\documentclass[11pt,leqno]{article}
\usepackage{amsfonts}
\usepackage{mathrsfs}
\usepackage[margin=1in]{geometry}
\usepackage{amssymb}
\usepackage{amsthm}
\usepackage{latexsym}
\usepackage{amsmath}
\usepackage{hyperref}
\newtheorem{thm}{Theorem}[subsection]
\newtheorem{defn}[thm]{Definition} 
\newtheorem{lem}[thm]{Lemma} 
\newtheorem{cor}[thm]{Corollary}
\newtheorem{pro}[thm]{Proposition}
\newtheorem{exa}[thm]{Example}

\allowdisplaybreaks
\newtheorem{rem}{\bf{Remark}}[subsection]
\newcommand{\s}{\sqrt}
\newcommand{\m}{\mathbb}
\newcommand{\R}{\mathbb{R}}
\newcommand{\N}{\mathbb{N}}

\begin{document}

\title{Long-time existence for semi-linear Klein-Gordon equations with quadratic potential}
\author{
Q.-D. Zhang \thanks{The author is supported by NSFC 10871175.}\ \thanks{email address: zjuzqd@163.com}\\
Department of Mathematics, Zhejiang University,
\\Hangzhou, 310027, China \\ \& \\ Universit{\'e} Paris 13, Institut Galil{\'e}e,\\
CNRS, UMR  7539, Laboratoire Analyse G{\'e}om{\'e}trie et Applications\\
99, Avenue J.-B. Cl{\'e}ment,\\
F-93430 Villetaneuse}
\date{}


\maketitle
\begin{abstract}
We prove that small smooth solutions of semi-linear Klein-Gordon equations with quadratic
potential exist over a longer interval than the one given by local existence theory, for
almost every value of mass. We use normal form for the Sobolev energy. The difficulty in
comparison with some similar results on the sphere comes from the fact that two
successive eigenvalues $\lambda, \lambda'$ of $\sqrt{-\Delta+|x|^2}$ may be separated by
a distance as small as $\frac{1}{\sqrt{\lambda}}$.
\end{abstract}


\addtocounter{section}{-1}
\section{Introduction}

Let $-\Delta + \vert x\vert^2$ be the harmonic oscillator on $\R^d$. This paper is
devoted to the proof of lower bounds for the existence time of solutions of non-linear
Klein-Gordon equations of type
\[ \begin{split}
(\partial_t^2 -\Delta + \vert x\vert^2 +m^2)v &= v^{\kappa+1}\\
v\vert_{t=0} &= \epsilon v_0\\
\partial_t v\vert_{t=0} &= \epsilon v_1
\end{split}\]
where $m\in \R^*_+$, $x^{\alpha}\partial_x^\beta v_j \in L^2$ when $\vert \alpha\vert +
\vert \beta\vert\leq s+1-j$ ($j=0,1)$ for a large enough integer $s$, and where
$\epsilon>0$ is small enough.

The similar equation without the quadratic potential $\vert x\vert^2$, and with data
small, smooth and compactly supported, has global solutions when $d\geq 2$ (see
Klainerman~\cite{K} and Shatah~\cite{Sh} for dimensions $d\geq 3$, Ozawa, Tsutaya and
Tsutsumi~\cite{OTT} when $d=2$). The situation is drastically different when we replace
$-\Delta$ by $-\Delta+\vert x\vert^2$, since the latter operator has pure point spectrum.
This prevents any time decay for solutions of the linear equation. Because of that, the
question of long time existence for  Klein-Gordon equations associated to the harmonic
oscillator is similar to the corresponding problem on compact manifolds.

For the equation $(\partial_t^2 -\Delta + m^2)v= v^{\kappa+1}$ on the circle
$\mathbb{S}^1$, it has been proved by Bourgain~\cite{Bo} and Bambusi~\cite{Ba}, that for
almost every $m>0$, the above equation has solutions defined on intervals of length
$c_N\epsilon^{-N}$ for any $N\in \N$, if the data are smooth and small enough (see also
the lectures of Gr{\'e}bert~\cite{G}). These results have been extended to the sphere
$\mathbb{S}^d$ instead of $\mathbb{S}^1$ by Bambusi, Delort, Gr{\'e}bert and
Szeftel~\cite{BDGS}. A key property in the proofs  is the structure of the spectrum of
$\sqrt{-\Delta}$ on $\mathbb{S}^d$. It is made of the integers, up to a small
perturbation, so that the gap between two successive eigenvalues is bounded from below by
a fixed constant.

A natural question is to examine which lower bounds on the time of existence of solutions
might be obtained when the eigenvalues of the operator do not satisfy such a gap
condition. The problem has been addressed for $(\partial_t^2 -\Delta + m^2)v =
v^{\kappa+1}$ on the torus $\mathbb{T}^d$ when $d\geq 2$ by Delort~\cite{D}. It has been
proved that for almost every $m>0$, the solution of such an equation exists over an
interval of time of length bounded from below by $c\epsilon^{-\kappa(1+2/d)}$ (up to a
logarithm) and has Sobolev norms of high index bounded on such an interval. Note that two
successive eigenvalues $\lambda, \lambda'$ of $\sqrt{-\Delta}$ on $\mathbb{T}^d$ might be
separated by an interval of length as small as $c/\lambda$. A natural question is then to
study the same problem for a model for which separation of eigenvalues is intermediate
between the cases of the sphere and of the torus. The harmonic oscillator provides such a
framework, as the distance between two successive eigenvalues $\lambda, \lambda'$ of
$\sqrt{-\Delta + \vert x\vert^2}$ is of order $1/\sqrt{\lambda}$. Our goal is to exploit
this to get for the corresponding Klein-Gordon equation a lower bound of the time of
existence of order $c\epsilon^{-4\kappa/3}$ when $d\geq 2$ (and a slightly better bound
if $d=1$).

Note that the estimate we get for the time of existence is explicit (given by the
exponent $-4\kappa/3$) and independent of the dimension $d$. This is in contrast with the
case of the torus, where the gain $2/d$ on the exponent brought by the method  goes to
zero as $d\to +\infty$. The point is that when the dimension increases, the multiplicity
of the eigenvalues of $-\Delta + \vert x\vert^2$ grows, while the spacing between
different eigenvalues remains essentially the same.

The method we use is based, as for similar problems on the sphere and the torus, on
normal form methods. Such an idea has been introduced in the study of non-linear
Klein-Gordon equations on $\R^d$ by Shatah~\cite{Sh}, and is at the root of the results
obtained on $\mathbb{S}^1, \mathbb{S}^d, \mathbb{T}^d$ in~\cite{Bo, Ba, BG, BDGS, D}. In
particular, we do not need to use any KAM results, unlike in the study of periodic or
quasi-periodic solutions of semi-linear wave or Klein-Gordon equations. For such a line
of studies, we refer to the books of Kuksin~\cite{Ku1, Ku2} and Craig~\cite{Cr} in the
case of the equation on $\mathbb{S}^1$, to Berti and Bolle~\cite{BB} for recent results
on the sphere, and to Bourgain~\cite{Bo2} and Elliasson-Kuksin~\cite{EK} in the case of
the torus.

Finally let us mention that very recently Gr{\'e}bert, Imekraz and Paturel~\cite{GP} have
studied the non-linear Schr{\"o}dinger equation associated to the harmonic oscillator.
They have obtained almost global existence of small solutions for this equation.

\numberwithin{equation}{subsection}
\section{The semi-linear Klein-Gordon equation}
\subsection{Sobolev Spaces} We introduce in this subsection Sobolev spaces we will work
with. From now on, we denote by $P=\sqrt{-\Delta+|x|^2}, x\in \mathbb{R}^d, d\ge 1$. The
operator $P^2=-\Delta+|x|^2$ is called the harmonic oscillator on $\mathbb{R}^d$. The
eigenvalues of $P^2$ are given by $\lambda_n^2$, where
\begin{equation}\label{eqn-B1}\lambda_n=\sqrt{2n+d}, n\in \mathbb{N}.\end{equation}
Let $\Pi_n$ be the orthogonal projector to the eigenspace associated to $\lambda_n^2$.
There are several ways to characterize these spaces. Of course we will show they are
equivalent after giving definitions.

\begin{defn}\label{def-B1} Let $s\in\mathbb{R}$. We define $\mathscr{H}_1^s(\mathbb{R}^d)$
to be the set of all functions $u\in L^2(\mathbb{R}^d)$ such that \\$
(\lambda_n^{s}||\Pi_nu||_{L^2})_{n\in\mathbb{N}}\in \ell^2$, equipped with the norm
defined by $||u||^2_{\mathscr{H}_1^s}=\sum_{n\in\mathbb{N}}\lambda_n^{2s}||
\Pi_nu||_{L^2}^2$.  \end{defn} The space $\mathscr{H}_1^s(\mathbb{R}^d)$ is the domain of
the operator $g(P)$ on $L^2(\mathbb{R}^d)$, which is defined using functional calculus
and  where
\begin{equation}\label{eqn-B2} g(r)=(1+r^2)^{\frac s2}, r\in \mathbb{R}.\end{equation}
Because of (\ref{eqn-B1}), we have
\begin{equation}\label{eqn-B3} ||g(P)u||_{L^2}\sim ||u||_{\mathscr{H}_1^s}.\end{equation}
\begin{defn}\label{def-B2} Let $s\in \mathbb{N}$. We define $\mathscr{H}_2^s(\mathbb{R}^d)$
to be the set of all functions $u\in L^2(\mathbb{R}^d)$ such that $x^\alpha\partial^\beta
u\in L^2 (\mathbb{R}^d), \forall   |\alpha|+|\beta|\le s$, equipped with the norm defined
by $||u||^2_{\mathscr{H}_2^s}=\sum_{|\alpha|+|\beta|\le s}||x^\alpha\partial^\beta
u||^2_{L^2}$.\end{defn} We shall give another definition of the space in the view point
of pseudo-differential theory. Let us first list some results from \cite{He}.
\begin{defn}\label{def-B3}We denote by $\Gamma^s(\mathbb{R}^{d})$, where $s\in
\mathbb{R}$, the set of all functions $u\in C^{\infty}(\mathbb{R}^{d})$ such that:
$\forall \alpha\in \mathbb{N}^{d}$, $\exists\  C_\alpha$, s.t.  $ \forall z\in
\mathbb{R}^{d}$, we have $|\partial_z^\alpha a(z)|\le C_\alpha \langle z
\rangle^{s-|\alpha|}$, where $\langle z\rangle=(1+|z|^2)^{\frac12}$ .\end{defn}

\begin{defn}\label{def-B4} Assume $a_j\in\Gamma^{s_j}(\mathbb{R}^{d})(j\in\mathbb{N}^*)$
and that $s_j$ is a decreasing sequence tending to $-\infty$. We say a function $a\in
C^\infty(\mathbb{R}^{d})$ satisfies: \begin{equation*}a\sim \sum_{j=1}^\infty
a_j\end{equation*} if: $\qquad \forall r\ge 2, r\in \mathbb{N}, \qquad
a-\sum_{j=1}^{r-1}a_j\in \Gamma^{s_r}(\mathbb{R}^{d})$.\end{defn}

We now would like to consider operators of the form
\begin{equation}\label{eqn-B4}Au(x)=(2\pi)^{-d}\iint e^{i(x-y)\cdot \xi}a(x, \xi)u(y)dyd\xi
\end{equation}
where $a(x, \xi)\in \Gamma^s(\mathbb{R}^{2d})$. We can also consider a more general
formula for the action of the operator
\begin{equation}\label{eqn-B5}Au(x)=(2\pi)^{-d}\iint e^{i(x-y)\cdot \xi}a(x, y,
\xi)u(y)dyd\xi\end{equation} where the function $a(x,y, \xi)$ is called the amplitude. We
will describe the class of amplitudes as following:
\begin{defn}\label{def-B5} Let $s\in \mathbb{R}$ and \ $\Omega^s(\mathbb{R}^{3d})$
denote the set of functions $a(x, y, \xi)\in C^\infty(\mathbb{R}^{3d})$, which for some
$s'\in \mathbb{R}$ satisfy
\begin{equation*} |\partial^\alpha_\xi\partial^\beta_x\partial^\gamma_y a(x,y,\xi)|\le
C_{\alpha\beta\gamma}\langle z\rangle^{s-(|\alpha|+|\beta|+|\gamma|)}\langle
x-y\rangle^{s'+|\alpha|+|\beta|+|\gamma|},\end{equation*} where $z=(x, y, \xi)\in
\mathbb{R}^{3d}.$\end{defn} The following proposition is a special case of proposition
1.1.4 in \cite{He}.
\begin{pro}\label{thm-B1}If $b\in \Gamma^s(\mathbb{R}^{2d})$, then $a( x, y, \xi)=b( x, \xi)$
and $a( x, y, \xi)=b( y, \xi)$ belong to
$\Omega^s(\mathbb{R}^{3d})$.\end{pro} Let $\chi(x, y, \xi)\in
C_0^\infty(\mathbb{R}^{3d}), \chi(0,0,0)=1$. It is shown by lemma 1.2.1 in \cite{He} that
(\ref{eqn-B5}) makes sense in the following way:
\begin{equation}\label{eqn-B6} Au(x)=\lim_{\varepsilon\rightarrow+0}(2\pi)^{-d}\iint
e^{i(x-y)\cdot\xi}\chi(\varepsilon x, \varepsilon y,
\varepsilon \xi)a(x, y, \xi)u(y)dyd\xi\end{equation}  if $a(x, y, \xi)\in
\Omega^s(\mathbb{R}^{3d})$ for some $s$. It is also shown in the same section of it the
operator $A$ is continuous from $\mathcal {S}(\mathbb{R}^d)$ to $\mathcal
{S}(\mathbb{R}^d)$ and it can be uniquely extended to an operator from $\mathcal
{S}'(\mathbb{R}^d)$ to $\mathcal {S}'(\mathbb{R}^d)$.

\begin{defn}\label{def-B6} The class of pseudo-differential operators $A$ of the form
(\ref{eqn-B5}) with amplitudes
$a\in \Omega^s(\mathbb{R}^{3d})$ will be denoted by $G^s(\mathbb{R}^d)$.\end{defn} We set
$G^{-\infty}(\mathbb{R}^d)=\bigcap_{s\in \mathbb{R}}G^s(\mathbb{R}^d)$.
\begin{exa} For $s\in \mathbb{N}$, the constant coefficient differential operator
$\sum\limits_{|\alpha|+|\beta|\le
s}c_{\alpha\beta}x^\alpha\partial^\beta$ is in the class $G^s(\mathbb{R}^d)$. \end{exa}
The class $G^s(\mathbb{R}^d)$ has some properties which are just theorems 1.3.1, 1.4.7,
1.4.8 in \cite{He}:
\begin{thm}\label{thm-B2} Let $s_1, s_2 \in \mathbb{R}$ and $A\in G^{s_1}(\mathbb{R}^d)$,
$A'\in G^{s_2}(\mathbb{R}^d)$. Then
$A\circ A'\in G^{s_1+s_2}(\mathbb{R}^d)$.\end{thm}
\begin{thm}\label{thm-B3} The operator $A\in G^0(\mathbb{R}^d)$ can be extended to a
bounded operator on $L^2(\mathbb{R}^d)$.\end{thm}
\begin{thm}\label{thm-B4} The operator $A\in G^s(\mathbb{R}^d)$ for $s<0$ can be extended
to a compact operator on
$L^2(\mathbb{R}^d)$. \end{thm} We shall give a subclass of that of pseudo-differential operators.
\begin{defn}\label{def-B7}We say $a\in \Gamma^{s}_{cl}(\mathbb{R}^d)$ if $a\in \Gamma^s
(\mathbb{R}^d)$ and $a$ has asymptotic expansion:
\begin{equation*}a\sim \sum_{j\in
\mathbb{N}}a_{s-j} \end{equation*} with $a_{s-j}\in C^\infty(\mathbb{R}^d)$ satisfying for
$\theta\ge 1, |x|+|\xi|\ge 1$
$$a_{s-j}(\theta x, \theta \xi)=\theta^{s-j}a_{s-j}(x, \xi).$$  \end{defn}

\begin{defn}\label{def-B8} Let $A$ be a pseudo-differential operator with amplitude $a\in
\Gamma^{s}_{cl}(\mathbb{R}^d)$. We then call $a_s$ defined above the principle symbol of
$A$.\end{defn}

\begin{defn}\label{def-B9} We say a pseudo-differential operator $A\in G^{s}_{cl}(\mathbb{R}^d)$
if its amplitude $a\in \Gamma^{s}_{cl}(\mathbb{R}^{2d})$. \end{defn} By proposition
\ref{thm-B1}, definition \ref{def-B9} is meaningful.

\begin{defn}\label{def-B10}We say that $A\in G^{s}_{cl}(\mathbb{R}^d)$ is globally
elliptic if we have: $\exists R>0, \exists\ C>0$ such that $\forall (x,\xi)\in
\mathbb{R}^{2d}$ satisfying $ |x|+|\xi|\ge R$, we have $|a_s(x, \xi)|\ge C(|x|+|\xi|)^s$,
where $a_s$ denotes the principle symbol of $A$.\end{defn} We can invert the operator
$A\in G^s_{cl}(\mathbb{R}^d)$ up to a regularizing operator, which is just theorem 1.5.7
in \cite{He}.
\begin{thm}\label{thm-B5} Let $A\in G^s_{cl}(\mathbb{R}^d)$ be a globally elliptic operator.
Then there is an operator $B\in G^{-s}_{cl}(\mathbb{R}^d)$ such that
\begin{equation}\label{eqn-B7}B\circ A=I+R_1, \qquad A\circ B=I+R_2,\end{equation}
where $R_1, R_2$ are regularizing, i.e. $R_1, R_2\in
G^{-\infty}(\mathbb{R}^d)$. \end{thm}

\begin{defn}\label{def-B11} Let $A$ be a pseudo-differential operator whose symbol is
$\langle \xi,x\rangle^s$ modulo $\Gamma^{s-1}_{cl}$.
 We define $\mathscr{H}_3^s(\mathbb{R}^d)$ to be the set of all functions $u\in \mathcal {S}'
 (\mathbb{R}^d)$ such that $Au\in
L^2(\mathbb{R}^d)$, equipped with the norm defined by $||u||^2_{\mathscr{H}_3^s}=
||Au||^2_{L^2}+||u||^2_{L^2}$. \end{defn}
\begin{rem}The pseudo-differential operator $A$ defined above is globally elliptic.
Thus by theorem \ref{thm-B5} if $Au\in L^2(\mathbb{R}^d)
$, we must have $u\in L^2(\mathbb{R}^d)$.\end{rem}
\begin{rem}$\mathscr{H}_3^s(\mathbb{R}^d)$ does not depend on the choice of $A$ according
to corollary 1.6.5 in \cite{He}. \end{rem}

\begin{cor} When $s\in \mathbb{N}$, definitions \ref{def-B1}, \ref{def-B2} and \ref{def-B11}
characterize the same space. Moreover
$\mathscr{H}_3^s(\mathbb{R}^d)=\mathscr{H}_1^s(\mathbb{R}^d)$ for any $s\in\mathbb{R}$.\end{cor}
\begin{proof} First let $s\in\mathbb{N}$. Since $A$ in definition \ref{def-B11} is globally elliptic,
by theorem \ref{thm-B5} there is $B\in
G^{-s}_{cl}(\mathbb{R}^d)$ such that
\begin{equation}\label{eqn-B8}B\circ A=I +R_1,\qquad A\circ B=I+R_2\end{equation} where $R_1, R_2$
are regularizing. Thus for any
$\alpha, \beta$ with $|\alpha|+|\beta|\le s$, by the example after definition
\ref{def-B7} and theorems \ref{thm-B2}, \ref{thm-B3} and \ref{thm-B4}, we have
$||x^\alpha\partial^\beta u||_{L^2}\le ||x^\alpha\partial^\beta
BAu||_{L^2}+||x^\alpha\partial^\beta R_1u||_{L^2}\le C(||Au||_{L^2}+||u||_{L^2})$, which
implies $||u||_{\mathscr{H}_2^s}\le C||u||_{\mathscr{H}_3^s}$. The inverse inequality
follows from the proof of proposition 1.6.6 in \cite{He}. Let us now prove that
definition \ref{def-B1} is equivalent to definition \ref{def-B11} for any
$s\in\mathbb{R}$.

By Theorem 1.11.2 in \cite{He} the operator $g(P)$ defined in (\ref{eqn-B2}) is an
essentially self-adjoint globally elliptic operator in the class $G^s(\mathbb{R}^d)$. We
have again by theorem \ref{thm-B5} that there is $Q\in G^{-s}_{cl}(\mathbb{R}^d)$ such
that
\begin{equation}\label{eqn-B9}g(P)\circ Q=I +R'_1,\qquad Q\circ g(P)=I+R'_2\end{equation}
where $R'_1, R'_2$ are regularizing.
We compute using (\ref{eqn-B3}), (\ref{eqn-B8}), (\ref{eqn-B9}) together with theorem
\ref{thm-B2} and theorem \ref{thm-B3}
\begin{equation*}\begin{split}||u||_{\mathscr{H}_1^s}\sim||g(P)u||_{L^2}&\le ||(g(P)\circ
B\circ A)u||_{L^2}+||(g(P)\circ R_1)u||_{L^2}
\\&\le C(||Au||_{L^2}+||u||_{L^2}) \le C||u||_{\mathscr{H}_3^s}
\end{split}\end{equation*} and \begin{equation*}\begin{split}||u||_{\mathscr{H}_3^s}&\le
C(||(A\circ Q\circ g(P))u||_{L^2}+||(A\circ R'_2)u||_{L^2}+||u||_{L^2})\\&\le
C(||g(P)u||_{L^2}+||u||_{L^2})\le C||u||_{\mathscr{H}_1^s},
\end{split}\end{equation*} where the last inequality follows from the fact $\lambda_n\ge 1$. \end{proof}

We denote $\mathscr{H}^s(\mathbb{R}^d)=\mathscr{H}_1^s(\mathbb{R}^d)=\mathscr{H}_3^s
(\mathbb{R}^d)$ when $s\in\mathbb{R}$. When
$s\in \mathbb{N}$, this space coincides with $\mathscr{H}_2^s(\mathbb{R}^d)$.
Let us present some properties of the spaces we
shall use.

\begin{pro}\label{thm-B7} If $s_1\le s_2$, then $\mathscr{H}^{s_2}(\mathbb{R}^d)
\hookrightarrow \mathscr{H}^{s_1}(\mathbb{R}^d).$\end{pro}
\begin{pro}\label{thm-B8} If $s>d/2$, then $\mathscr{H}^s(\mathbb{R}^d) \hookrightarrow
 L^\infty(\mathbb{R}^d) $.\end{pro}
\begin{pro}\label{thm-B9} Let $f\in C^\infty(\mathbb{R}), f(0)=0, u\in \mathscr{H}^s(\mathbb{R}^d),
 s\in \mathbb{N}, s>d$. Then we have
$f(u)\in \mathscr{H}^s(\mathbb{R}^d)$. Moreover if $f$ vanishes at order $p+1$ at 0,
where $p\in \mathbb{N}$, then $||f(u)||_{\mathscr{H}^s}\le C||u||^{p+1}_{\mathscr{H}^s}$.
\end{pro}
\begin{proof} Proposition \ref{thm-B7} and \ref{thm-B8} follow respectively from the definition
and Sobolev embedding. By the chain rule, for $|\alpha|+|\beta|\le s$,
$x^\alpha\partial^\beta f(u)$ may be written as the sum of terms of following form:
\begin{equation*} x^\alpha f^{(k)}(u)(\partial^{\beta_1}u)\dots(\partial^{\beta_k}u),
\end{equation*} where $k\le s, |\alpha|+\sum_{i=1}^{k}|\beta_i|\le s, |\beta_i|>0, i=1,
\dots, k$. Let $j_0$ be the index such that $|\beta_{j_0}|$ is the largest among
$|\beta_1|, \dots, |\beta_k|$. Thus we must have $|\beta_i|\le \frac{s}{2}$, $i\neq
{j_0}$. By the assumption on $s$ and proposition \ref{thm-B8}, $\partial^\gamma u\in
L^\infty(\mathbb{R}^d)$ if $|\gamma|\le \frac{d}{2}$. We then estimate the factor
$x^\alpha\partial^{\beta_{j_0}}u$ of the above quantities in $L^2$-norm and others in
$L^\infty$-norm. Thus we have $f(u)\in {\mathscr{H}}^s(\mathbb{R}^d)$ by proposition
\ref{thm-B8}. When $f$ vanishes at 0 at order $p+1$, by Taylor formula there is a smooth
function $h$ such that $f(u)=u^{p+1}h(u)$. Then we argue as above to get an upper bound
of $||f(u)||_{\mathscr{H}^s}$ by $C||u||_{\mathscr{H}^s}^{p}||u||_{\mathscr{H}^s}$. This
concludes the proof.
\end{proof}
\begin{rem}\label{thm-B6} Proposition \ref{thm-B9} actually holds true for $s>d/2$ if
we argue as the proof of corollary 6.4.4 in \cite{Ho}. Since we will consider only in
$\mathscr{H}^s(\mathbb{R}^d)$ for large $s$, the lower bound of $s$ is not important.
\end{rem}

\numberwithin{equation}{subsection}
\subsection{Statement of main theorem}\label{sec}
Let $d$ be an integer, $d\ge 1$ and $F: \mathbb{R}\rightarrow \mathbb{R}$ a real valued
smooth function vanishing at order $\kappa+1$ at $0$, $\kappa\in \mathbb{N}^*$. Let $m\in
\R^*_+$. we consider the solution $v$ of the following Cauchy problem:
     \begin{eqnarray} \label{eqn-C1}\left\{\begin{split}
     (\partial_t^2-\Delta+|x|^2+m^2)v&=F(v) \hspace{3mm}\mbox{on}\hspace{2mm}[-T,T]\times\mathbb{R}^d\\
                                          v(0,x)&=\epsilon v_0\\
                               \partial_t v(0,x)&=\epsilon v_1,
     \end{split}\right.\end{eqnarray}
where $v_0\in \mathscr{H}^{s+1}(\mathbb{R}^d), v_1\in \mathscr{H}^s(\mathbb{R}^d)$, and
$\epsilon >0$ is a small parameter. By local existence theory one knows  that if $s$ is
large enough and $\epsilon \in (0,1)$, equation (\ref{eqn-C1}) admits for any $(v_0,
v_1)$ in the unit ball of $\mathscr{H}^{s+1}(\mathbb{R}^d)\times
\mathscr{H}^s(\mathbb{R}^d)$ a unique smooth solution defined on the interval $|t|\le
c\epsilon^{-\kappa}$, for some uniform positive constant c. Moreover, $||v(t,\cdot
)||_{\mathscr{H}^{s+1}}+||\partial_tv(t,\cdot )||_{\mathscr{H}^s}$ may be controlled by
$C\epsilon$, for another uniform constant $C>0$, on the interval of existence. The goal
would be to obtain existence over an interval of longer length under convenient condition
by controlling the Sobolev energy. Our main result is the following:

\begin{thm}\label{2-1-1}
There is a zero measure subset $\mathcal {N}$ of $\R^*_+$ and for every $m\in
\R^*_+-\mathcal {N}$, there are $\epsilon_0>0, c>0, s_0\in\mathbb{N}$ such that for any
$s\ge s_0, s\in \mathbb{N}, \epsilon\in(0,\epsilon_0)$, any pair $(v_0, v_1)$ of real
valued functions belonging to the unit ball of
$\mathscr{H}^{s+1}(\mathbb{R}^d)\times\mathscr{H}^{s}(\mathbb{R}^d)$, problem
(\ref{eqn-C1}) has a unique solution
\begin{equation} u\in C^0((-T_\epsilon, T_\epsilon), \mathscr{H}^{s+1}(\mathbb{R}^d))\cap C^1
((-T_\epsilon, T_\epsilon),
\mathscr{H}^s(\mathbb{R}^d)),\end{equation} where $T_\epsilon$ has a lower bound
$T_\epsilon\ge c\epsilon^{-\frac{4}{3}(1-\rho)\kappa}$ for any $\rho>0$ if $d\ge2$ and
$T_\epsilon\ge c\epsilon^{-\frac{25}{18}(1-\rho)\kappa}$ for any $\rho>0$ if $d=1$.
Moreover, the solution is uniformly bounded in $\mathscr{H}^{s+1} (\mathbb{R}^d)$ on
$(-T_\epsilon, T_\epsilon)$ and $\partial_t u$ is uniformly bounded in
$\mathscr{H}^s(\mathbb{R}^d)$ on the same interval.
  \end{thm}

\numberwithin{equation}{subsection}
\subsection{ A property of spectral projectors on $\mathbb{R}^d$}
\paragraph{}
As we have pointed out $P$ has eigenvalues given by $\lambda_n=\sqrt{2n+d}, n\in
\mathbb{N}$. Remark that $\Pi_n$ is the orthogonal projector of $L^2(\mathbb{R}^d)$ onto
the eigenspace associated to $\lambda_n^2$. Let us first introduce some notations. For
$\xi_0, \xi_1, \dots, \xi_{p+1}$ $p+2$ nonnegative real numbers, let $\xi_{i_0},
\xi_{i_1}, \xi_{i_2}$ be respectively the largest, the second largest and the third
largest elements among them and $\xi'$ the largest element among $\xi_1, \dots, \xi_{p}$,
that is,
\begin{equation}\begin{split} \label{eqn-D1}\xi_{i_0}=\max\{\xi_0, \dots, \xi_{p+1}\},
\qquad \xi_{i_1}=\max(\{\xi_0, \dots,
\xi_{p+1}\}-\{\xi_{i_0}\}), \\ \xi_{i_2}=\max(\{\xi_1, \dots, \xi_{p+1}\}-\{\xi_{i_0},
\xi_{i_1}\} ),\qquad \xi'=\max\{\xi_1,
\dots, \xi_p\}. \end{split}\end{equation} Denote
\begin{equation}\label{eqn-D2}\mu(\xi_0, \dots, \xi_{p+1})=(1+\sqrt{\xi_{i_1}})(1+\sqrt{\xi_{i_2}}).
\end{equation} Set also
\begin{equation}\label{eqn-D3}S(\xi_0, \dots, \xi_{p+1})=|\xi_{i_0}-\xi_{i_1}|+\mu(\xi_0, \dots,
\xi_{p+1}).\end{equation}

The main result of this subsection is the following one:
\begin{thm}\label{thm-D1} There is a  $\nu\in \R^*_+$, depending only on $p\ (p\in \mathbb{N}^*)$
and dimension $d$, and for any $N\in\mathbb{N}$, there is a $C_N>0$ such that for any
$n_0, \dots, n_{p+1}   \in \mathbb{N}$, any $u_0,\dots,u_{p+1} \in L^2(\mathbb{R}^d)$,
\begin{equation}\label{eqn-D4}|\int \Pi_{n_0}u_0\dots \Pi_{n_{p+1}}u_{p+1}dx|\le C_N
(1+\sqrt{n_{i_2}})^\nu\frac{\mu(n_0, \dots, n_{p+1})^{N}}{S(n_0,
\dots, n_{p+1})^N}\prod_{j=0}^{p+1}||u_j||_{L^2}.\end{equation} Furthermore if $d=1$, we
may find for any $\varsigma\in (0,1)$
\begin{equation}\label{eqn-D5}|\int \Pi_{n_0}u_0\dots \Pi_{n_{p+1}}u_{p+1}dx|\le
C_N\frac{(1+\sqrt{n_{i_2}})^\nu}{(1+\sqrt{n_{i_0}})^{\frac{1}{6}(1-\varsigma)}}\frac{\mu(n_0,
\dots, n_{p+1})^{N}}{S(n_0, \dots,
n_{p+1})^N}\prod_{j=0}^{p+1}||u_j||_{L^2}.\end{equation}  \end{thm}

\begin{proof} By the symmetries we may assume $n_0\ge n_1\ge\dots\ge n_{p+1}$. Then recalling
the definition of $\lambda_{n}$ in
(\ref{eqn-B1}), we only need to show under the condition of theorem \ref{thm-D1}
\begin{equation}\label{eqn-D6}|\int
\Pi_{n_0}u_0\dots \Pi_{n_{p+1}}u_{p+1}dx|\le
C_N\lambda_{n_2}^\nu\frac{(\lambda_{n_1}\lambda_{n_2})^{N}}{(|\lambda_{n_0}^2-\lambda_{n_1}^2|
+\lambda_{n_1}\lambda_{n_2})^N}\prod_{j=0}^
{p+1}||u_j||_{L^2}\end{equation} and when $d=1$
\begin{equation}\label{eqn-D7}|\int \Pi_{n_0}u_0\dots \Pi_{n_{p+1}}u_{p+1}dx|\le C_N
\frac{\lambda_{n_2}^{\nu}}{\lambda_{n_0}^{\frac{1}{6}
(1-\varsigma)}}\frac{(\lambda_{n_1}\lambda_{n_2})^{N}}{(|\lambda_{n_0}^2-\lambda_{n_1}^2|
+\lambda_{n_1}\lambda_{n_2})^N}
\prod_{j=0}^{p+1}||u_j||_{L^2}\end{equation} for any $\varsigma\in (0,1)$. We follow the
proof of proposition 3.6 in \cite{GP}. Let $A$ be a linear operator which maps
$D(P^{2k})$ into itself. We define a sequence of operators
\begin{equation}\label{eqn-D8} A_N=[P^2, A_{N-1}]; \qquad A_0=A.\end{equation} Then using
integration by parts we have
\begin{equation}\label{eqn-D9}(\lambda_{n_0}^2-\lambda_{n_1}^2)^N\langle A\ \Pi_{n_1}u_1,
\Pi_{n_0}u_0\rangle =\langle
A_N\Pi_{n_1}u_1, \Pi_{n_0}u_0\rangle. \end{equation} Now we set $A$ to be the
multiplication operator generated by the function
$$a(x)=(\Pi_{n_2}u_2)\dots(\Pi_{n_{p+1}}u_{p+1}).$$ Then an induction argument shows
\begin{equation}\label{eqn-D10} A_N=\sum_{|\beta|+|\gamma|\le N,\
|\alpha|+|\beta|+|\gamma|\le 2N}C_{\alpha\beta\gamma}(\partial^\alpha
a)x^\beta\partial^\gamma\end{equation} for constants $C_{\alpha\beta\gamma}$. Therefore
we compute for some $\nu\thinspace'>\frac{d}{2}$
\begin{equation}\label{eqn-D11}\begin{split}&\qquad |(\lambda_{n_0}^2-\lambda_{n_1}^2)^N\int
(\Pi_{n_0}u_0)\dots(\Pi_{n_{p+1}}u_{p+1})dx|\\& \le C\sum_{|\beta|+|\gamma|\le N,\
|\alpha|+|\beta|+|\gamma|\le 2N}||(\partial^\alpha
a)x^\beta\partial^\gamma\Pi_{n_1}u_1||_{L^2}||\Pi_{n_0}u_0||_{L^2}\\ &\le
C\sum_{|\beta|+|\gamma|\le N,\ |\alpha|+|\beta|+|\gamma|\le
2N}||a||_{\mathscr{H}^{\nu\thinspace'+|\alpha|}}||\Pi_{n_1}u_1||_{\mathscr{H}^{|\beta|+|\gamma|}}
||\Pi_{n_0}u_0||_{L^2},
\end{split}\end{equation} where in the last estimate we used definition \ref{def-B2}
and proposition \ref{thm-B8}. Remark that by definition
\ref{def-B1}, one has for any $s\ge 0$
\begin{equation}\label{eqn-D12}||\Pi_{n}u||_{\mathscr{H}^{s}}\le
C\lambda_{n}^s||\Pi_{n}u||_{L^2}.\end{equation} This estimate together with the proof of
proposition \ref{thm-B9} gives for $n_2\ge n_3 \dots \ge n_{p+1}$
\begin{equation}\label{eqn-D13} ||a||_{\mathscr{H}^{\nu'+|\alpha|}}\le
C\lambda^{\nu+|\alpha|}_{n_{2}}\prod_{j=2}^{p+1}||\Pi_{n_j}u_j||_{L^2}\end{equation} for
some $\nu>0$ depending only on $p$ and dimension $d$. Thus we have
\begin{equation}\label{eqn-D14}\begin{split}&\qquad
|(\lambda_{n_0}^2-\lambda_{n_1}^2)^N\int (\Pi_{n_0}u_0)\dots(\Pi_{n_{p+1}}u_{p+1})dx|\\&
\le C\sum_{|\beta|+|\gamma|\le N,\ |\alpha|+|\beta|+|\gamma|\le
2N}\lambda^{\nu+|\alpha|}_{n_{2}}\lambda_{n_1}^{|\beta|+|\gamma|}\prod_{j=0}^{p+1}||\Pi_{n_j}u_j||
_{L^2}\\&\le C\sum_{|\alpha|\le
N}\lambda^{\nu+2N-|\alpha|}_{n_{2}}\lambda_{n_1}^{|\alpha|}\prod_{j=0}^{p+1}||\Pi_{n_j}u_j||_{L^2}
\\ & \le
C\lambda_{n_{2}}^{\nu+2N}(\frac{\lambda_{n_1}}{\lambda_{n_2}})^N\prod_{j=0}^{p+1}||\Pi_{n_j}u_j||_{L^2}\\
& \le
C\lambda_{n_2}^{\nu}(\lambda_{n_1}\lambda_{n_2})^N\prod_{j=0}^{p+1}||\Pi_{n_j}u_j||_{L^2}.
\end{split}\end{equation} Now if $\lambda_{n_1}\lambda_{n_2}\le |\lambda_{n_0}^2-\lambda_{n_1}^2|$,
then the last estimate implies (\ref{eqn-D6}), while if $\lambda_{n_1}\lambda_{n_2}>
|\lambda_{n_0}^2-\lambda_{n_1}^2|$, then
$\frac{\lambda_{n_1}\lambda_{n_2}}{|\lambda_{n_0}^2-\lambda_{n_1}^2|+\lambda_{n_1}\lambda_{n_2}}\ge\frac12$
and thus (\ref{eqn-D6}) is trivially true.

On the other hand, we use the property of the eigenfunctions (see\cite{KT}), which in
dimension $d=1$ says that if $\phi_n$ is the eigenfunction associated to $\lambda_n^2$,
then one has $||\phi_n||_{L^\infty}\le C\lambda_n^{-\frac{1}{6}}$. Therefore we have
\begin{equation}\label{eqn-D15}||\Pi_nu||_{L^\infty}\le
C\lambda_{n}^{-\frac{1}{6}}||\Pi_nu||_{L^2}\end{equation} since in this case the
eigenvalues are simple. This estimate gives us \begin{equation}\label{eqn-D16}|\int
\Pi_{n_0}u_0\dots \Pi_{n_{p+1}}u_{p+1}dx|\le
C\lambda_{n_0}^{-\frac{1}{6}}\prod_{j=0}^{p+1}||\Pi_{n_j}u_j||_{L^2}.\end{equation}
Combining (\ref{eqn-D16}) with (\ref{eqn-D6}) one gets (\ref{eqn-D7}) for all $N\ge1$ and
some $\nu>0$ in the case $d=1$. This concludes the proof. \end{proof}

\numberwithin{equation}{subsection}
\section{Long time existence}
\subsection{Definition and properties of multilinear operators}

Denote by $\mathcal {E}$ the algebraic direct sum of the ranges of the
$\Pi_n\thinspace{'}s,n\in\mathbb{N}$. With notations (\ref{eqn-D1}), (\ref{eqn-D2}) and
(\ref{eqn-D3}) we give the following definition.
\begin{defn}\label{def-E1}
Let $\nu\in\R_+,\ \tau \in \m{R}, \ p\in\mathbb{N}^*$. We denote by $\mathcal
{M}_{p+1}^{\nu,\tau}$ the space of all $p+1$-linear operators
$(u_1,\dots,u_{p+1})\rightarrow M(u_1,\dots,u_{p+1})$, defined on $\mathcal {E}\times
\dots\times\mathcal {E}$ with values in $L^2(\mathbb{R}^d)$ such that
\newcounter{Lcount}{\usecounter{Lcount}}
\begin{itemize}
\item  For every $(n_0,\dots,n_{p+1})\in\mathbb{N}^{p+2},u_1,\dots,u_{p+1}\in\mathcal{E}$
       \begin{equation}\label{eqn-E1}\Pi_{n_0}[M(\Pi_{n_1}u_1,\dots,\Pi_{n_{p+1}}u_{p+1})]=0, \end{equation}
       if $|n_0-n_{p+1}|>\frac12(n_0+n_{p+1})$ or $n'\overset{def}{=}\max\{n_1, \dots, n_p\}>n_{p+1}$.
\item  For any $N\in\mathbb{N}$, there is a $C>0$ such that for every $(n_0,\dots,n_{p+1})\in\mathbb{N}^{p+2}$,\\
$u_1,\dots,u_{p+1}\in\mathcal{E}$, one has
\begin{equation}\label{eqn-E2}\begin{split}
&||\Pi_{n_0}[M(\Pi_{n_1}u_1,\dots,\Pi_{n_{p+1}}u_{p+1})]||_{L^2}\\ &\qquad\qquad\le
C(1+\s{n_0}+\s{n_{p+1}})^{\tau}(1+\sqrt{n'})^\nu\frac{\mu(n_0, \dots, n_{p+1})^N}{S(n_0,
\dots, n_{p+1})^N}\prod_{j=1}^{p+1}||u_j||_{L^2}. \end{split}\end{equation}
\end{itemize}
The best constant in the preceding inequality will be denoted by $||M||_{{\mathcal {M}_{p+1,N}^{\nu,\tau}}}$.
\end{defn}


We may extend the operators in $\mathcal {M}_{p+1}^{\nu,\tau}$ to Sobolev spaces.

\begin{pro}\label{thm-E1} Let $\nu\in\R_+,\  \tau\in \mathbb{R},\ p\in \mathbb{N}^*,\ s\in\mathbb{N},\ s>\nu+3$.
Then any element  $M\in\mathcal {M}_{p+1}^{\nu,\tau}$ extends as a bounded operator from
$\mathscr{H}^s(\mathbb{R}^d)\times\dots\times \mathscr{H}^s(\mathbb{R}^d)$ to $
\mathscr{H}^{s-\tau-1}(\mathbb{R}^d)$. Moreover, for any $s_0\in (\nu+3, s]$, there is
$C>0$ such that for any $M\in\mathcal {M}_{p+1}^{\nu,\tau}$, and any $u_1, \dots,
u_{p+1}\in \mathscr{H}^s(\mathbb{R}^d)$,
\begin{equation}\label{eqn-E3}
||M(u_1, \dots, u_{p+1})||_{\mathscr{H}^{s-\tau-1}}\leq C||M||_{\mathcal
{M}_{p+1,N}^{\nu,\tau}}\sum_{j=1}^{p+1}\biggr[||u_j||_{\mathscr{H}^s}\prod_{k\neq
j}||u_k||_{\mathscr{H}^{s_0}}\biggr].
\end{equation}
\end{pro}

\begin{proof}  The proof is a modification of proposition 4.4 in \cite{DS1}.
There is one derivative lost compared to that case. We give it for the convenience of the
reader. Using definition \ref{def-B1} we write
\begin{multline}\label{eqn-E4}||M(u_1, \dots, u_{p+1})||_{\mathscr{H}^{s-\tau-1}}^2\\ \le
     C\sum_{n_0}||\sum_{n_1}\dots\sum_{n_{p+1}}\Pi_{n_0}M(\Pi_{n_1}u_1,\dots,\Pi_{n_{p+1}}u_{p+1})||_{L^2}^2
     (1+\sqrt{n_0})^{2s-2\tau-2}
\end{multline}
Because of (\ref{eqn-E1}) and using the symmetries we may assume
\begin{equation}\label{eqn-E5} n_0\sim n_{p+1} \qquad \text{and}\qquad n_1\le\dots\le
n_p\le n_{p+1}\le Cn_0 \end{equation} when estimating the above quantity. Consequently,
we have
\begin{equation}\label{eqn-E6}\begin{split} \mu(n_0, \dots, n_{p+1})\sim (1+\sqrt{n_p})(1+\sqrt{n_{p+1}}),\\ S(n_0, \dots,
n_{p+1})\sim |n_0-n_{p+1}|+\mu(n_0, \dots, n_{p+1}).\end{split} \end{equation} By
(\ref{eqn-E2}) the square root of the general term over $n_0$ sum in (\ref{eqn-E4}) is
smaller than
\begin{equation}\label{eqn-E7}
     C\sum_{n_1\le \dots \le n_{p+1}}(1+\sqrt{n_0})^{s-1}\frac{(1+\sqrt{n_p})^{\nu}\mu(n_0,
      \dots, n_{p+1})^N}{S(n_0, \dots, n_{p+1})
     ^N}\prod_1^{p+1}||\Pi_{n_j}u_j||_{L^2}.
\end{equation}
We have by (\ref{eqn-E5}) and (\ref{eqn-E6})
\begin{equation}\label{eqn-E7'}\frac{\mu(n_0, \dots, n_{p+1})}{S(n_0, \dots, n_{p+1})}
\sim \frac{1+\sqrt{n_p}}{|\sqrt{n_0}-\sqrt{n_{p+1}}|+1+\sqrt{n_p}}.\end{equation} The
following fact will be useful in this section: For $q\in\mathbb{N}, A\ge 1$ and $N>1$,
there is a $C>0$ independent of $q$ and $A$ such that
\begin{equation}\label{eqn-E8}
          \sum_{n\in\mathbb{N}}\frac{1}{(|\sqrt{n}-\sqrt{q}|+A)^N}\le C\frac{1+\sqrt{q}}{A^{N-2}}.
\end{equation}
Let $\iota>2$ be a constant as close to 2 as wanted. Using (\ref{eqn-E7'}) and
(\ref{eqn-E8}) we deduce
         \begin{equation}\label{eqn-E9}\begin{split}
            \sum_{n_0}\frac{\mu(n_0, \dots, n_{p+1})^{\iota}}{S(n_0, \dots, n_{p+1})^{\iota}}
            \le C(1+\sqrt{n_{p+1}})(1+\sqrt{n_p})^2,\\ \sum_{n_{p+1}}\frac{\mu(n_0, \dots,
            n_{p+1})^{\iota}}{S(n_0, \dots, n_{p+1})^{\iota}}
            \le C(1+\sqrt{n_{0}})(1+\sqrt{n_p})^2.
    \end{split}   \end{equation}
We estimate the sum over $n_1\le\dots\le n_{p+1}$ in (\ref{eqn-E7}) by
\begin{equation}\begin{split}\label{eqn-E10} &\qquad C\biggr(\sum_{{n_1}\le\dots\le
     n_{p+1}}\frac{(1+\sqrt{n_p})^\nu\mu^{\iota}}{S^\iota}\prod_{j=1}^p||\Pi_{n_j}u_j||_{L^2}\biggr)^{1/2}\\ & \times
     \biggr(\sum_{{n_1}\le\dots\le
    n_{p+1}}(1+\sqrt{n_0})^{2s-2}(1+\sqrt{n_p})^\nu\frac{\mu^{2N-\iota}}{S^{2N-\iota}}\prod_{j=1}^p||\Pi_{n_j}u_j||_{L^2}
    ||\Pi_{n_{p+1}}u_{p+1}||_{L^2}^2 \biggr)^{1/2}.
\end{split} \end{equation}
Using (\ref{eqn-E9}) to handle $n_{p+1}$ sum, we bound the first factor in
(\ref{eqn-E10}) from above by
$C(1+\sqrt{n_0})^{\frac12}\Pi_{j=1}^p||u_j||^{\frac12}_{\mathscr{H}^{s_0}}$ if
$s_0>\nu+3$ using definition \ref{def-B1}. Incorporating $(1+\sqrt{n_0})^{\frac12}$ into
the second factor, we have to bound the quantity
\begin{equation}\label{c10}
\biggr(\sum_{{n_1}\le\dots\le
n_{p+1}}(1+\sqrt{n_0})^{2s-1}(1+\sqrt{n_p})^\nu\frac{\mu^{2N-\iota}}{S^{2N-\iota}}\prod_{j=1}^p||
\Pi_{n_j}u_j||_{L^2}||\Pi_{n_{p+1}}u_{p+1}||_{L^2}^2\biggr)^{1/2}.
\end{equation}   By (\ref{eqn-E5}) and $\mu\le S$
we have
\begin{equation}\label{eqn-E11}(1+\sqrt{n_0})^{2s-1}(\frac{\mu}{S})^{2N-\iota}\le C(1+\sqrt{n_{p+1}})^
{2s-1}(\frac{\mu}{S})^\iota\end{equation}
if $N>\iota$. Plugging in  (\ref{c10}), (\ref{eqn-E10}) and then (\ref{eqn-E4}) we bound
from above the $n_0$ sum in (\ref{eqn-E4}) by
\begin{equation}\label{eqn-E12}   C\prod_1^p||u_j||_{\mathscr{H}^{s_0}}\sum_{n_1\le \dots \le n_{p+1}\le
Cn_0}(1+\sqrt{n_{p+1}})^{2s-1}(1+\sqrt{n_p})^\nu (\frac{\mu}{S})^\iota\prod_{j=1}^p
||\Pi_{n_j}u_j||_{L^2}||\Pi_{n_{p+1}}u_{p+1}||_{L^2}^2.\end{equation} Changing the order
of sums for $n_0$ and $n_{p+1}$, we then use (\ref{eqn-E9}) to handle $n_0$ sum and get a
control of (\ref{eqn-E12}) by
$C\prod_{j=1}^p||u_j||^2_{\mathscr{H}^{s_0}}||u_{p+1}||^2_{\mathscr{H}^s}$ according to
definition \ref{def-B1} if $s
>\nu+3$. This concludes the proof.
\end{proof}


\paragraph{}

Let us define convenient subspaces of the spaces of definition \ref{def-E1}.
\begin{defn}\label{def-E2}
Let $\nu\in  \R_+,\ \tau\in\mathbb{R},\ p\in \mathbb{N}^*,\omega : \{0, \dots, p+1\}
\rightarrow \{-1, 1\}$ be given.
\begin{itemize}
\item  If $\sum_{j=0}^{p+1}\omega(j)\neq 0$, we set $\widetilde{\mathcal{M}}_{p+1}^{\nu,\tau}(\omega)=\mathcal{M}_
                    {p+1}^{\nu,\tau};$
\item  If $\sum_{j=0}^{p+1}\omega(j)= 0$, we denote by $\widetilde{\mathcal{M}}_{p+1}^{\nu,\tau}(\omega)$
the closed subspace of
$\mathcal{M}_{p+1}^{\nu,\tau}$ given by those $M\in \mathcal{M}_{p+1}^{\nu,\tau}$ such
that
\begin{equation}\label{eqn-E13}\Pi_{n_0}M(\Pi_{n_1}u_1, \dots, \Pi_{n_{p+1}}u_{p+1})\equiv 0 \end{equation}
for any $(n_0, \dots, n_{p+1})\in \m{N}^{p+2}$ such that there is a bijection $\sigma$
from $\{j; 0\le j \le p+1, \omega(j)=-1\}$ to $\{j;0\le j\le p+1,\omega(j)=1\}$ so that
for any $j$ in the first set $n_{\sigma{(j)}}=n_j$.
\end{itemize}
\end{defn}

We shall have to use also classes of remainder operators. If $n_1, \dots, n_{p+1}\in \N$
and $j_0\in\{1,\dots, p+1\}$ is such that $n_{j_0}=\max\{n_1, \dots, n_{p+1}\}$, we
denote
\begin{equation}\label{eqn-E14'} \max\nolimits_2(\sqrt{n_1}, \dots, \sqrt{n_{p+1}})=
1+\max\{\sqrt{n_j}\thinspace ;\ 1\le j \le p+1, j\neq j_0\}.\end{equation}

\begin{defn}\label{def-E3}
Let $ \nu \in \R_+,\ \tau\in \mathbb{R},\ p\in \mathbb{N}^*$. We denote by
$\mathcal{R}_{p+1}^{\nu, \tau}$ the space of $\mathbb{C}\  (p+1)$-linear maps from
$\mathcal {E}\times\dots\times \mathcal {E}\rightarrow L^2(\m{R}^d), (u_1, \dots,
u_{p+1})\rightarrow R(u_1, \dots, u_{p+1})$ such that for any $N\in \m{N}$, there is a
$C>0$ such that for any $(n_0, \dots, n_{p+1}) \in \m{N}^{p+2},$ any $u_1, \dots, u_{p+1}
\in \mathcal {E}$,
\begin{equation}\label{eqn-E14}
        ||\Pi_{n_0}R(\Pi_{n_1}u_1, \dots, \Pi_{n_{p+1}}u_{p+1})||_{L^2}\le
C(1+\sqrt{n_0})^\tau\frac{\max\nolimits_2(\sqrt{n_1}, \dots, \sqrt{n_{p+1}})^{\nu+N}}
{(1+\s{n_0}+\dots+\sqrt{n_{p+1}})^{N}}
        \prod_{j=1}^{p+1}||u_j||_{L^2}.
       \end{equation}
\end{defn}

\paragraph{}
The elements in $\mathcal {R}_{p+1}^{\nu, \tau}$ also extend as bounded operators on
Sobolev spaces.

\begin{pro}\label{thm-E2} Let $\nu \in \R_+,\ \tau\in \mathbb{R},\ p\in \mathbb{N}^*$ be given.
There is $s_0\in \m{N}$ such that for any $s\ge s_0$, any $R\in \mathcal
{R}_{p+1}^{\nu,\tau}$, $(u_1, \dots, u_{p+1})\rightarrow R(u_1, \dots, u_{p+1})$ extends
as a bounded map from $\mathscr{H}^s(\m{R}^d)\times \dots \times
\mathscr{H}^s(\m{R}^d)\rightarrow \mathscr{H}^{2s-\nu-\tau-7}(\mathbb{R}^d)$. Moreover
one has
\begin{equation}\label{eqn-E15}
     ||R(u_1, \dots, u_{p+1})||_{\mathscr{H}^{2s-\nu-\tau-7}}\le C\sum_{1\le j_1< j_2\le p+1}
     \biggr[||u_{j_1}||_{\mathscr{H}^s}||u_{j_2}||_{\mathscr{H}^s}\prod_{k\neq
    j_1, k\neq j_2}||u_k||_{\mathscr{H}^{s_0}}\biggr].
\end{equation}
\end{pro}

\begin{proof} We may assume $\tau=0$. By definition \ref{def-B1} we have to bound
$||\Pi_{n_0}R(u_1, \dots, u_{p+1})||_{L^2}$ from above by
$(1+\sqrt{n_0})^{-2s+\nu+7}c_{n_0}$ for a sequence $(c_{n_0})_{n_0}$ in $\ell^2$. To do
that we decompose $u_j$ as $\sum_{n_j}\Pi_{n_j}u_j$ and use (\ref{eqn-E14}). By symmetry
we limit ourselves to summation over
\begin{equation}\label{eqn-E16}n_1\le \dots\le n_{p+1},\end{equation}
from which we deduce \begin{equation}\label{eqn-E17} \max\nolimits_2(\sqrt{n_1}, \dots,
\sqrt{n_{p+1}})=1+\sqrt{n_p} .\end{equation} Therefore we are done if we can bound from
above
\begin{equation}\label{eqn-E18}
   C\sum_{n_1\le\dots\le n_{p+1}}\frac{(1+\s{n_p})^{\nu+N}}{(1+\s{n_0}+\dots+\s{n_{p+1}})^{N}}
   \prod_{j=1}^{p-1}(1+\s{n_j})^{-s_0}(1+\s{n_{p}})^{-s}(1+\s{n_{p+1}})^{-s}
\end{equation} by $(1+\sqrt{n_0})^{-2s+\nu+7}c_{n_0}$ for $s_0, s$
large enough with $s\ge s_0$ since $||\Pi_{n_j}u_j||_{L^2}
\le C(1+\sqrt{n_j})^{-s}||u_j||_{\mathscr {H}^s}$. Using (\ref{eqn-E16}) we get an upper
bound of (\ref{eqn-E18}) by
\begin{equation}\label{eqn-E19}C\sum_{n_1\le\dots\le n_{p+1}}\frac{(1+\s{n_p})^{\nu+N-2s}}{(1+\s{n_0}+\s{n_{p+1}})^{N}}
\prod_{j=1}^{p-1}(1+\s{n_j})^{-s_0}\end{equation} Using the fact $\sum_{n\in
\mathbb{N}}\frac{1}{(\sqrt{n}+A)^N}\le \frac{C}{A^{N-2}}$ for $N>2$ and $A\ge 1$ , we
take the sum over $n_{p+1}$ to get an upper bound of (\ref{eqn-E18}) by
\begin{equation}\label{eqn-E20}C\sum_{n_1\le\dots\le n_{p}}\frac{(1+\s{n_p})^{\nu+N-2s}}{(1+\s{n_0})^{N-2}}
\prod_{j=1}^{p-1}(1+\s{n_j})^{-s_0}\end{equation} if $N>2$. Now take $N=2s-\nu-\frac52$
and sum over $n_1, \dots, n_p$. This gives the upper bound we want and thus concludes the
proof. \end{proof}

\begin{defn}\label{def-E4} Let $\nu\in  \R_+,\ \tau\in\mathbb{R},\ p\in \mathbb{N}^*,\omega : \{0, \dots,
p+1\} \rightarrow \{-1, 1\}$ be given.
\begin{itemize}
\item  If $\sum_{j=0}^{p+1}\omega(j)\neq 0$, we set $\widetilde{\mathcal{R}}_{p+1}^{\nu,\tau}(\omega)=\mathcal{R}_
                    {p+1}^{\nu,\tau};$
\item  If $\sum_{j=0}^{p+1}\omega(j)= 0$, we denote by $\widetilde{\mathcal{R}}_{p+1}^{\nu,\tau}(\omega)$
the closed subspace of
$\mathcal{R}_{p+1}^{\nu,\tau}$ given by those $R\in \mathcal{M}_{p+1}^{\nu,\tau}$ such
that
\begin{equation}\label{eqn-E21}\Pi_{n_0}R(\Pi_{n_1}u_1, \dots, \Pi_{n_{p+1}}u_{p+1})\equiv 0 \end{equation}
for any $(n_0, \dots, n_{p+1})\in \m{N}^{p+2}$ such that there is a bijection $\sigma$
from $\{j; 0\le j \le p+1, \omega(j)=-1\}$ to $\{j;0\le j\le p+1,\omega(j)=1\}$ so that
for any $j$ in the first set $n_{\sigma{(j)}}=n_j$.
\end{itemize}
\end{defn}



\subsection{Rewriting of the equation and the energy}
In this subsection we will write the time derivative of the energy in terms of
multilinear operators defined in the previous subsection. To do that, we shall need to
analyze the nonlinearity. Decompose \begin{equation}\label{eqn-F1}
-F(v)=-\sum_{p=\kappa}^{2\kappa-1}\frac{\partial^{p+1}_vF(0)}{(p+1)!}v^{p+1}+G(v)\end{equation}
where $G(v)$ vanishes at order $2\kappa+1$ at $v=0$. One has
$$cv^{p+1}=c\sum_{n_1}\dots\sum_{n_{p+1}}(\Pi_{n_1}v)\dots(\Pi_{n_{p+1}}v)$$ for a real
constant $c$. One may also write this as $A_p(v)\cdot v$ where $A_p(v)$ is an operator of
form
\begin{equation}\label{s}A_p(v)\cdot w=\sum_{n_1}\dots\sum_{n_{p+1}}B(n_1, \dots,
n_{p+1})(\Pi_{n_1}v)\dots(\Pi_{n_p}v)(\Pi_{n_{p+1}}w),\end{equation} where $B(n_1, \dots,
n_{p+1})$ is a real valued bounded function supported on $\max\{n_1,\dots,n_p\}\le
n_{p+1}$ and $B$ is constant valued on the domain $\max\{n_1,\dots,n_p\}< n_{p+1}$. For
instance, when $p=2$, one may write
$$\{(n_1,n_2,n_3);n_j\in \mathbb{N}\}=\{\max\{n_1,n_2\}\le n_3\}\cup\{n_1\ge n_2\ \text{and}
\ n_1>n_3\}\cup\{n_1<n_2\ \text{and}\
n_2>n_3\}$$ and
\begin{equation*}\begin{split}\sum_{n_1}&\sum_{n_2}\sum_{n_3}(\Pi_{n_1}v)(\Pi_{n_2}v)(\Pi_{n_3}v)
=\sum\textbf{1}_{\{\max\{n_ 1,n_2\}\le n_3\}}(\Pi_{n_1}v)
(\Pi_{n_2}v)(\Pi_{n_3}v)\\&+\sum\textbf{1}_{\{n_3\ge n_2\ \text{and}\
n_3>n_1\}}(\Pi_{n_1}v)(\Pi_{n_2}v)(\Pi_{n_3}v)+\sum \textbf{1}_{\{n_3> n_2\ \text{and}\
n_3>n_1\}}(\Pi_{n_1}v)(\Pi_{n_2}v)(\Pi_{n_3}v)\end{split}\end{equation*} using the
symmetries, so that in this case
$$B(n_1, n_2, n_3)=c(\textbf{1}_{\{\max\{n_1,n_2\}\le n_3\}}+\textbf{1}_{\{n_3\ge n_2
\ \text{and}\ n_3>n_1\}}+\textbf{1}_ {\{n_3> n_2\ \text{and}\ n_3>n_1\}}).$$ So if we
make a change of unknown $u=(D_t+\Lambda_m)v$ with $$D_t=-i\partial_t, \qquad
\Lambda_m=\sqrt{-\Delta+|x|^2+m^2},$$ we may write using (\ref{eqn-F1})
\begin{equation}\label{eqn-F2}(D_t-\Lambda_m)u=-\sum_{p=\kappa}^{2\kappa-1}A_p\biggr(
\Lambda_m^{-1}(\frac{u+\bar{u}}{2})\biggr)
\Lambda_m^{-1}(\frac{u+\bar{u}}{2})+G\biggr(\Lambda_m^{-1}(\frac{u+\bar{u}}{2})\biggr).\end{equation}
Denote
$C(u,\bar{u})=-\frac12\sum^{2\kappa-1}_{p=\kappa}A_p\biggr(\Lambda_m^{-1}(\frac{u+\bar{u}}{2})\biggr)\Lambda_m^{-1}$
so that
\begin{equation}\label{eqn-F3}(D_t-\Lambda_m)u=C(u,\bar{u})u+C(u,\bar{u})\bar{u}+G\biggr(\Lambda_m^{-1}
(\frac{u+\bar{u}}{2})\biggr).\end{equation}
We have to estimate for the solution $u$ of (\ref{eqn-F2})\begin{equation}\label{eqn-F4}
\Theta_s(u(t,\cdot))=\frac12\langle\Lambda_m^s u(t, \cdot), \Lambda_m^s u(t, \cdot)\rangle.\end{equation}

Now comes the main result of this subsection:
\begin{pro}\label{thm-F1} There are $\nu\in \R_+$ and large enough $s_0$ such that
for any natural number $s\ge s_0$, there are: \begin{itemize}
\item Multilinear operators $M^p_\ell\in \widetilde{\mathcal {M}}^{\nu, 2s-a}_{p+1}(\omega_\ell)$, $\kappa\le p\le
2\kappa-1$, $ 0\le \ell\le p$ with $\omega_\ell$ defined by $\omega_\ell(j)=-1$, $j=0,
\dots, \ell,\ \omega_\ell(j)=1, j=\ell+1, \dots, p+1$ and $a=2$ if $d\ge2$ and
$a=\frac{13}{6}-\varsigma$ for any $\varsigma\in (0,1)$ if $d=1$;
\item Multilinear operators $\widetilde{M}^p_\ell\in \widetilde{\mathcal {M}}^{\nu, 2s-1}_{p+1}
(\widetilde{\omega}_\ell)$,
$\kappa\le p\le 2\kappa-1$, $ 0\le \ell\le p$ with $\widetilde{\omega}_\ell$ defined by
$\widetilde{\omega}_\ell(j)=-1$, $j=0, \dots, \ell, p+1$,$\ \widetilde{\omega}_\ell(j)=1,
j=\ell+1, \dots, p$;
\item Multilinear operators $R^p_\ell\in\widetilde{\mathcal {R}}^{\nu, 2s}_{p+1}(\omega_\ell),
\widetilde{R}^p_\ell\in \widetilde{\mathcal {R}}^{\nu,
2s}_{p+1}(\widetilde{\omega}_\ell)$, $\kappa\le p\le 2\kappa-1$, $ 0\le \ell\le p$;
\item A map $u\rightarrow T(u)$ defined on $\mathscr{H}^s(\mathbb{R}^d)$ with values in $\mathbb{R}$,
satisfying when $||u||_{\mathscr {H}^s}\le 1$, $|T(u)|\le
C||u||^{2\kappa+2}_{\mathscr{H}^s}$
\end{itemize}
such that \begin{equation}\label{eqn-F5}\begin{split}
&\frac{d}{dt}\Theta_s(u(t,\cdot))=\sum_{p=\kappa} ^{2\kappa-1}\sum_{\ell=0}^{p}Re\
i\langle M^p_\ell(\underset{\ell}{\underbrace{\bar{u}, \dots, \bar{u}}},
\underset{p+1-\ell}{\underbrace{u, \dots, u}}), u\rangle\\&+\sum_{p=\kappa}
^{2\kappa-1}\sum_{\ell=0}^{p} Re\ i\langle
\widetilde{M}^p_\ell(\underset{\ell}{\underbrace{\bar{u}, \dots, \bar{u}}},
\underset{p-\ell}{\underbrace{u, \dots, u}},\bar{u}), u\rangle+\sum_{p=\kappa}
^{2\kappa-1}\sum_{\ell=0}^{p}Re\ i\langle R^p_\ell(\underset{\ell}{\underbrace{\bar{u},
\dots, \bar{u}}}, \underset{p+1-\ell}{\underbrace{u, \dots, u}}), u\rangle\\&+
\sum_{p=\kappa} ^{2\kappa-1}\sum_{\ell=0}^{p}Re\ i\langle
\widetilde{R}^p_\ell(\underset{\ell}{\underbrace{\bar{u}, \dots, \bar{u}}},
\underset{p-\ell}{\underbrace{u, \dots, u}}, \bar{u}),
u\rangle+T(u).\end{split}\end{equation}
\end{pro}

\begin{proof} We compute according to (\ref{eqn-F3}) \begin{multline}\label{eqn-F6}\frac{d}{dt}\Theta_s
(u(t, \cdot))=Re\ i\langle \Lambda_m^sD_tu, \Lambda_m^s u\rangle\\=Re\ i\langle
\Lambda_m^sC(u,\bar{u})u, \Lambda_m^s u\rangle + Re\ i\langle \Lambda_m^s
C(u,\bar{u})\bar{u}, \Lambda_m^s u\rangle +Re\ i\langle
\Lambda_m^sG(\Lambda_m^{-1}(\frac{u+\bar{u}}{2})), \Lambda_m^{s}u\rangle.
\end{multline} The last term in the right hand side of (\ref{eqn-F6}) contributes to the last term in (\ref{eqn-F5})
by proposition \ref{thm-B9}. Let us treat the other two terms in the right hand side of
(\ref{eqn-F6}).

\begin{lem}\label{thm-F2} There are $M^p_\ell\in
\widetilde{\mathcal {M}}^{\nu, 2s-a}_{p+1}(\omega_\ell)$, $R^p_\ell\in\widetilde{\mathcal
{R}}^{\nu, 2s}_{p+1}(\omega_\ell)$, $\kappa\le p\le 2\kappa-1$, $ 0\le \ell\le p$ with
$\omega_\ell$ defined by $\omega_\ell(j)=-1$, $j=0, \dots, \ell,\ \omega_\ell(j)=1,
j=\ell+1, \dots, p+1$ and $a=2$ if $d\ge2$ and $a=\frac{13}{6}-\varsigma$ for any
$\varsigma\in (0,1)$ if $d=1$, such that
\begin{equation}\begin{split}\label{eqn-F7}Re\ i\langle \Lambda_m^sC(u,\bar{u})u,
\Lambda_m^s u\rangle=\sum_{p=\kappa} ^{2\kappa-1}\sum_{\ell=0}^{p}Re\ i\langle
M^p_\ell(\underset{\ell}{\underbrace{\bar{u}, \dots, \bar{u}}},
\underset{p+1-\ell}{\underbrace{u, \dots, u}}), u\rangle\\+\sum_{p=\kappa}
^{2\kappa-1}\sum_{\ell=0}^{p}Re\ i\langle R^p_\ell(\underset{\ell}{\underbrace{\bar{u},
\dots, \bar{u}}}, \underset{p+1-\ell}{\underbrace{u, \dots, u}}), u\rangle.\end{split}
\end{equation} \end{lem}

\begin{proof}[Proof of Lemma \ref{thm-F2}:]   Let $\chi$ be a cut-off
function near $0$ with small support and $\lambda_n$ defined in (\ref{eqn-B1}). We may
decompose the operator $A_p(v)$ defined in (\ref{s}) as \begin{equation}\label{eqn-F21}
A_p(v)=A_p^1(v)+A_p^2(v)+A_p^3(v),\end{equation} where $A_p^j(v)(j=1,2,3)$ are operators
of form \begin{equation}\label{eqn-F22}\begin{split} A_p^1(v)\cdot
w=\sum_{n_0}\dots\sum_{n_{p+1}}B_1(n_0,\dots,n_{p+1})\Pi_{n_0}[(\Pi_{n_1}v)\dots(\Pi_{n_p}v)(\Pi_{n_{p+1}}w)],\\
A_p^2(v)\cdot
w=\sum_{n_0}\dots\sum_{n_{p+1}}B_2(n_0,\dots,n_{p+1})\Pi_{n_0}[(\Pi_{n_1}v)\dots(\Pi_{n_p}v)(\Pi_{n_{p+1}}w)],\\
A_p^3(v)\cdot
w=\sum_{n_1}\dots\sum_{n_{p+1}}B_3(n_1,\dots,n_{p+1})\Pi_{n_0}[(\Pi_{n_1}v)\dots(\Pi_{n_p}v)(\Pi_{n_{p+1}}w)],
\end{split}\end{equation} with \begin{equation}\label{eqn-F23}\begin{split}&B_1(n_0,\dots,n_{p+1})=B(n_1,\dots,n_{p+1})
\chi\biggr(\frac{|\lambda^2_{n_0}-\lambda_{n_{p+1}}^2|}{\lambda^2_{n_0}+\lambda_{n_{p+1}}^2}\biggr)
\textbf{1}_{\{\max\{n_1,\dots,n_p\}<\delta n_{p+1}\}},\\&
B_2(n_0,\dots,n_{p+1})=B(n_1,\dots,n_{p+1})
\biggr(1-\chi\biggr(\frac{|\lambda^2_{n_0}-\lambda_{n_{p+1}}^2|}{\lambda^2_{n_0}+\lambda_{n_{p+1}}^2}\biggr)\biggr)
\textbf{1}_{\{\max\{n_1,\dots,n_p\}<\delta n_{p+1}\}},\\&
B_3(n_1,\dots,n_{p+1})=B(n_1,\dots,n_{p+1})\textbf{1}_{\{\max\{n_1,\dots,n_p\}\ge\delta
n_{p+1}\}},\end{split}\end{equation} with some small $\delta>0$.  Therefore for the
operator $C(u,\bar{u})$ defined above (\ref{eqn-F3}), we have
\begin{equation}\label{eqn-F24}C(u,\bar{u})=-\frac12\sum_{j=1}^3\sum_{p=\kappa}^{2\kappa-1}A_p^j\biggr(\Lambda_{m}^{-1}
(\frac{u+\bar{u}}{2})\biggr)\Lambda_m^{-1}.\end{equation} So the left hand side of
(\ref{eqn-F7}) may be written as
\begin{equation}\label{eqn-F25}-\frac12\sum_{j=1}^3\sum_{p=\kappa}^{2\kappa-1}Re\ i\langle
\Lambda_m^{2s}A_p^j\biggr(\Lambda_{m}^{-1}(\frac{u+\bar{u}}{2})\biggr)\Lambda_m^{-1}u,
u\rangle:=\sum_{j=1}^3\sum_{p=\kappa}^{2\kappa-1}I^j_p.\end{equation}

Let us treat these quantities term by term. \paragraph{}\textbf{(i)} The term $I^1_p$.

Note that $-4I^1_p$ equals to
\begin{equation}\label{eqn-F26}Re\ i\langle
\biggr[\Lambda_m^{2s}A^1_p\biggr(\Lambda_{m}^{-1}(\frac{u+\bar{u}}{2})\biggr)\Lambda_m^{-1}-
\biggr(A_p^1\biggr(\Lambda_{m}^{-1}(\frac{u+\bar{u}}{2})\biggr)\Lambda_m^{-1}\biggr)^*\Lambda_m^{2s}\biggr]u,u
\rangle,\end{equation} which may be written as
\begin{equation}\begin{split}\label{eqn-F27} &Re\ i\langle \biggr[\Lambda_m^{2s},
A^1_p\biggr(\Lambda_m^{-1}(\frac{u+\bar{u}}{2})\biggr)\Lambda_m^{-1}\biggr]u,
u\rangle\\&+Re\ i\langle
\biggr[A^1_p\biggr(\Lambda_m^{-1}(\frac{u+\bar{u}}{2})\biggr)\Lambda_m^{-1}-\biggr(A^1_p\biggr(\Lambda_m^{-1}
(\frac{u+\bar{u}}{2})\biggr)\Lambda_m^{-1}\biggr)^*\biggr]\Lambda_m^{2s}u, u\rangle:=I+II
\end{split}\end{equation} We expand the first term in (\ref{eqn-F27}) using (\ref{eqn-F22}) to get
\begin{align}\label{eqn-F28} I=Re\ i \langle \sum_{n\in\N^{p+2}}\pi_1
\Pi_{n_0}\biggr[\biggr(\Pi_{n_1}\Lambda_m^{-1}(\frac{u+\bar{u}}{2})\biggr)\dots\biggr(\Pi_{n_p}\Lambda_m^{-1}
(\frac{u+\bar{u}}{2})\biggr)\biggr(\Pi_{n_{p+1}}\Lambda_m^{-1}u\biggr)\biggr],u \rangle\\
\nonumber =Re\
i\langle\sum_{n\in\N^{p+2}}\sum_{\ell=0}^{p}\pi_2\Pi_{n_0}[(\Pi_{n_1}\Lambda_m^{-1}\bar{u})\dots
(\Pi_{n_\ell}\Lambda_m^{-1}\bar{u})
(\Pi_{n_{\ell+1}}\Lambda_m^{-1}u)\dots(\Pi_{n_{p+1}}\Lambda_m^{-1}u)],u
\rangle\\\nonumber=Re\ i\sum_{n\in\N^{p+2}}\sum_{\ell=0}^{p}\pi_2\int
(\Pi_{n_0}\bar{u})(\Pi_{n_1}\Lambda_m^{-1}\bar{u})\dots(\Pi_{n_\ell}\Lambda_m^{-1}\bar{u})
(\Pi_{n_{\ell+1}}\Lambda_m^{-1}u)\dots(\Pi_{n_{p+1}}\Lambda_m^{-1}u)dx,\end{align} where
we have used notations
\begin{gather}\label{eqn-F29}\nonumber n=(n_0,\dots,n_{p+1}),\\\pi_1=B_1(n_0,\dots,
n_{p+1})[(m^2+\lambda^2_{n_0})^s-(m^2+\lambda^2_{n_{p+1}})^s],\\\nonumber
\pi_2=\frac{1}{2^p}\begin{pmatrix} p\\ \ell \end{pmatrix} B_1(n_0,\dots,
n_{p+1})[(m^2+\lambda^2_{n_0})^s-(m^2+\lambda^2_{n_{p+1}})^s].\end{gather} Let
$\omega_\ell$ be defined in the statement of the lemma and set
\begin{align}\label{eqn-F30}\begin{split}
S_p^\ell=\{(n_0,\dots, n_{p+1})\in\N^{p+2}; \thinspace \mbox{ there exists a bijection }
\sigma  \mbox{ from } \\ \{j; 0\le j \le p+1, \omega_\ell(j)=-1\}  \mbox{ to } \{ j; 0\le
j\le p+1, \omega_\ell(j)=1\} \\ \mbox{ such that for each } j \mbox{ in the first set }
n_j=n_{\sigma(j)}\}.\end{split}\end{align} Now we look at the integral in the last line
of (\ref{eqn-F28}). If $n\in S^\ell_p$ with $S^\ell_p\ne\emptyset$, there is a bijection
$\sigma$ from $\{0,\dots, \ell\}$ to $\{\ell, \dots, p+1\}$ such that $n_j=n_{\sigma(j)},
j=0, \dots, \ell$. So we may couple $\Pi_{n_j}\bar{u}, j=0, \dots, \ell$ with
$\Pi_{n_{\sigma(j)}}u, j=0, \dots, \ell$. Since $\pi_2$ is real, we get zero if we take
the sum over $n\in S^\ell_p$ when computing the right hand side of (\ref{eqn-F28}).
Therefore we may assume $n\notin S^\ell_p$ when computing $I$. Now we define
\begin{equation}\label{eqn-F31}M^{p,1}_\ell(u_1,\dots, u_{p+1})=-\frac14\sum_{n\notin
S^\ell_p}\pi_2\Pi_{n_0}[(\Pi_{n_1}\Lambda_m^{-1}u_1)\dots(\Pi_{n_{p+1}}\Lambda_m^{-1}u_{p+1})].\end{equation}
It follows from the second equality in (\ref{eqn-F28}) that
\begin{equation}\label{eqn-F32}I=-4\sum_{\ell=0}^pRe\ i\langle
M^{p,1}_\ell(\underset{\ell}{\underbrace{\bar{u},\dots,\bar{u}}},\underset{p+1-\ell}{\underbrace{u,\dots,
u}}), u \rangle.\end{equation}  Let us turn to the term $II$ in (\ref{eqn-F27}). Note
that $A^1_p(v)^*$ is an operator of form
\begin{equation}\label{eqn-F33}A^1_p(v)^*\cdot w=\sum_{n\in \N^{p+2}}B_1(n_{p+1},n_1, \dots,n_p, n_{0})\Pi_{n_0}
[(\Pi_{n_1}v)\dots(\Pi_{n_p}v)(\Pi_{n_{p+1}}w)].\end{equation} Thus we may compute using
(\ref{eqn-F22})
\begin{equation}\label{eqn-F34}\begin{split}II=Re\ i\langle \sum_{n\in\N^{p+2}}\sum_{\ell=0}^p\pi_3\Pi_{n_0}
[(\Pi_{n_1}\Lambda_m^{-1}\bar{u})\dots(\Pi_{n_\ell}\Lambda_m^{-1}\bar{u})(\Pi_{n_{\ell+1}}\Lambda_m^{-1}u)\dots
(\Pi_{n_p}\Lambda_m^{-1}u)(\Pi_{n_{p+1}}\Lambda_m^{2s}u)], u \rangle\\= Re\
i\sum_{n\in\N^{p+2}}\sum_{\ell=0}^p\pi_3\int
(\Pi_{n_0}\bar{u})(\Pi_{n_1}\Lambda_m^{-1}\bar{u})\dots(\Pi_{n_\ell}\Lambda_m^{-1}\bar{u})(\Pi_{n_{\ell+1}}\Lambda_m^{-1}u)\dots
(\Pi_{n_p}\Lambda_m^{-1}u)(\Pi_{n_{p+1}}\Lambda_m^{2s}u)dx,\end{split}\end{equation}
where
\begin{equation}\label{eqn-F35}\pi_3=\frac{1}{2^p}\begin{pmatrix}p\\ \ell\end{pmatrix}[B_1(n_0, n_1,\dots, n_p,n_{p+1})
(m^2+\lambda^2_{n_{p+1}})^{-\frac12}-B_1(n_{p+1}, n_1,\dots,
n_p,n_{0})(m^2+\lambda_{n_0}^{2})^{-\frac12}].\end{equation} With the same reasoning as
in the paragraph above (\ref{eqn-F31}) we get zero if we take the sum over $n\in
S^\ell_p$ when computing the right hand side of (\ref{eqn-F34}). So we may assume
$n\notin S^\ell_p$ and define
\begin{equation}\label{eqn-F36}M^{p,2}_\ell(u_1,\dots,u_{p+1})=-\frac14\sum_{n\notin S^\ell_p}\pi_3\Pi_{n_0}
[(\Pi_{n_1}\Lambda_m^{-1}u_1)\dots(\Pi_{n_{p}}\Lambda_m^{-1}u_p)(\Pi_{n_{p+1}}\Lambda_m^{2s}u_{p+1})].\end{equation}
It follows from (\ref{eqn-F34}) that
\begin{equation}\label{eqn-F37}II=-4\sum_{\ell=0}^pRe\ i\langle
M^{p,2}_\ell(\underset{\ell}{\underbrace{\bar{u},\dots,\bar{u}}},\underset{p+1-\ell}{\underbrace{u,\dots,
u}}), u \rangle.\end{equation}

Let us check that $M^{p,1}_\ell, M^{p,2}_\ell\in\widetilde{\mathcal {M}}^{\nu,
2s-a}_{p+1}(\omega_\ell)$ for some $\nu>0$, where $a=2$ if $d\ge 2$ and
$a=\frac{13}{6}-\varsigma$ for any $\varsigma\in (0,1)$ if $d=1$. Since the function
$B_1(n_0,\dots,n_{n_{p+1}})$ is supported on domain $n'=\max\{n_1,\dots,n_p\}<\delta
n_{p+1}$ and $n_0\sim n_{p+1}$ (this is because of the cut-off function and
(\ref{eqn-B1})), we see that (\ref{eqn-E1}) holds true if $supp\chi$ and $\delta$ are
small. Let us use theorem \ref{thm-D1} to show that (\ref{eqn-E2}) holds true with
$\tau=2s-a$ for $M^{p,1}_\ell$ and $M^{p,2}_\ell$. Remark that we have
\begin{gather}\label{eqn-F38}
|\pi_2|\le C(1+|\sqrt{n_0}-\sqrt{n_{p+1}}|)(1+\sqrt{n_0}+\sqrt{n_{p+1}})^{2s-1},\\
\label{eqn-F39}|\pi_3| \le
C(1+\sqrt{n'})^2(1+|\sqrt{n_0}-\sqrt{n_{p+1}}|)(1+\sqrt{n_0}+\sqrt{n_{p+1}})^{-2}.\end{gather}
Indeed, (\ref{eqn-F38}) follows from the fact
$$|(m^2+\lambda^2_{n_0})^s-(m^2+\lambda^2_{n_{p+1}})^s|\le
C(|\lambda_{n_0}-\lambda_{n_{p+1}}|)(1+\lambda_{n_0}+\lambda_{n_{p+1}})^{2s-1}.$$ If
$n'<\delta n_0$ and $n'<\delta n_{p+1}$ for small $\delta>0$, then $$B_1(n_0,
n_1,\dots,n_p,n_{p+1})=B_1(n_{p+1},n_1,\dots,n_p,n_0)$$ since $B(n_1,\dots,n_{p+1})$ is
constant valued on the domain $n'<n_{p+1}$. Thus (\ref{eqn-F39}) follows from the fact
$$|(m^2+\lambda^2_{n_0})^{-\frac12}-(m^2+\lambda^2_{n_{p+1}})^{-\frac12}|\le
C(|\lambda_{n_0}-\lambda_{n_{p+1}}|)(1+\lambda_{n_0}+\lambda_{n_{p+1}})^{-2}.
$$ Otherwise, assume $n'\ge\delta n_0$ or $n'\ge\delta n_{p+1}$. Then we must have $n'\ge
Cn_0$ and $n'\ge Cn_{p+1}$ if $B_1$ is non zero, since $n_0\sim n_{p+1}$ which is because
of the cut-off function. In this case, (\ref{eqn-F39}) holds true trivially.

 Moreover, on the support of
$\Pi_{n_0}M^{p,l}_\ell(\Pi_{n_1}u_1,\dots,\Pi_{n_{p+1}}u_{p+1})(l=1,2)$, i.e., $n_0\sim
n_{p+1}$ and  $n_{p+1}\ge \max\{n_1,\dots, n_p\}=n'$, we have
\begin{gather}\label{eqn-F40}\nonumber 1+\sqrt{n_{i_2}}\sim 1+\sqrt{n'},\\
\mu(n_0, \dots, n_{p+1})\sim (1+\sqrt{n_{p+1}})(1+\sqrt{n'}),\\ \nonumber  S(n_0,\dots,
n_{p+1})\sim |n_0-n_{p+1}|+(1+\sqrt{n_{p+1}})(1+\sqrt{n'}), \end{gather} from which we
deduce
\begin{equation}\label{eqn-F41} \frac{\mu(n_0, \dots, n_{p+1})}{S(n_0,\dots,
n_{p+1})}\sim \frac{1+\sqrt{n'}}{|\sqrt{n_0}-\sqrt{n_{p+1}}|+1+\sqrt{n'}}.\end{equation}
Thus $$(1+|\sqrt{n_0}-\sqrt{n_{p+1}}|)\frac{\mu(n_0, \dots, n_{p+1})}{S(n_0,\dots,
n_{p+1})}\le C(1+\sqrt{n'}).$$ Then we use theorem \ref{thm-D1} (with dimension $d\ge2$)
to get for $l=1,2$
\begin{equation}\begin{split}\label{eqn-F42}&\qquad\qquad||\Pi_{n_0}M^{p,l}_\ell
(\Pi_{n_1}u_1,\dots,\Pi_{n_{p+1}}u_{p+1})||_{L^2}\\
&\le
C(1+\sqrt{n_0}+\sqrt{n_{p+1}})^{2s-2}(1+\sqrt{n'})^{\nu+2}(1+|\sqrt{n_0}-\sqrt{n_{p+1}}|)
\frac{\mu(n_0,\dots,n_{p+1})^N}{S(n_0, \dots,
n_{p+1})^N}\prod_{j=1}^{p+1}||u_j||_{L^2}\\
&\le
C(1+\sqrt{n_0}+\sqrt{n_{p+1}})^{2s-2}(1+\sqrt{n'})^{\nu+3}\frac{\mu(n_0,\dots,n_{p+1})^{N-1}}{S(n_0,
\dots, n_{p+1})^{N-1}}\prod_{j=1}^{p+1}||u_j||_{L^2}.\end{split}\end{equation} So
$M^{p,l}_\ell\in\mathcal {M}^{\nu,2s-2}_\ell$ for some other $\nu>0$ in dimension
$d\ge2$. The case of dimension one is similar. (\ref{eqn-E13}) with $\omega=\omega_\ell$
is satisfied by definition. Thus $M^{p,1}_\ell,M^{p,2}_\ell\in \widetilde{\mathcal
{M}}^{\nu, 2s-a}_{p+1}(\omega_\ell)$ and we have proved \begin{equation}
I^1_p=\sum_{\ell=0}^pRe\ i\langle
M^{p,1}_\ell(\underset{\ell}{\underbrace{\bar{u},\dots,\bar{u}}},\underset{p+1-\ell}{\underbrace{u,\dots,
u}}), u \rangle+\sum_{\ell=0}^pRe\ i\langle
M^{p,2}_\ell(\underset{\ell}{\underbrace{\bar{u},\dots,\bar{u}}},\underset{p+1-\ell}{\underbrace{u,\dots,
u}}), u \rangle.\end{equation}

\paragraph{}
\textbf{(ii)} The term $I^2_p$.

Using (\ref{eqn-F22}) we get
\begin{equation}\begin{split}\label{eqn-F44}&-2I^2_p=Re\ i\langle \sum_{n\in\N^{p+2}}\sum_{\ell=0}^p
\pi_4\Lambda_m^{2s}\Pi_{n_0}[(\Pi_{n_1}\Lambda_m^{-1}\bar{u})\dots
(\Pi_{n_\ell}\Lambda_m^{-1}\bar{u})(\Pi_{n_{\ell+1}}\Lambda_m^{-1}u)\dots(\Pi_{n_{p+1}}\Lambda_m^{-1}u)],u\rangle
\\&=Re\ i \sum_{n\in\N^{p+2}}\sum_{\ell=0}^p
\pi_4\int (\Pi_{n_0}\Lambda_m^{2s}\bar{u})(\Pi_{n_1}\Lambda_m^{-1}\bar{u})\dots
(\Pi_{n_\ell}\Lambda_m^{-1}\bar{u})(\Pi_{n_{\ell+1}}\Lambda_m^{-1}u)\dots(\Pi_{n_{p+1}}\Lambda_m^{-1}u)dx
\end{split}\end{equation} where $$\pi_4=\frac{1}{2^p}\begin{pmatrix}p\\
\ell\end{pmatrix}B_2(n_0,\dots,n_{p+1}).$$ We may rule out the sum over $n\in S^\ell_p$
in the above computation with the same reasoning as in the paragraph above
(\ref{eqn-F31}). Thus if we define
\begin{equation}\label{eqn-F45}R^{p,1}_\ell(u_1,\dots,u_{p+1})=-\frac12\sum_{n\notin S^\ell_p}\pi_4
\Lambda_m^{2s}\Pi_{n_0}
[(\Pi_{n_1}\Lambda_m^{-1}u_1)\dots(\Pi_{n_{p+1}}\Lambda_m^{-1}u_{p+1})],\end{equation} we
have \begin{equation}\label{eqn-F46}I^2_p=\sum_{\ell=0}^pRe\ i\langle R^{p,1}_\ell
(\underset{\ell}{\underbrace{\bar{u},\dots,\bar{u}}},\underset{p+1-\ell}{\underbrace
{u,\dots,u}}),u\rangle.\end{equation} From the support property of function $B_2(n_0,
\dots, n_{p+1})$ we know that \\
$\Pi_{n_0}R^{p,1}_\ell(\Pi_{n_1}u_1,\dots,\Pi_{n_{p+1}}u_{p+1})$ is supported on
$\max\{n_1,\dots,n_p\}<\delta n_{p+1}$ and $|n_0-n_{p+1}|\ge c(n_0+n_{p+1})$ for some
small $c>0$. Therefore, on its support,
 if $n_0>C n_{p+1}$ for
a large $C$, we have
\begin{gather*}\mu(n_0, \dots,n_{p+1})=(1+\sqrt{n_{p+1}})(1+\sqrt{n'})\le (1+\sqrt{n_{0}})(1+\sqrt{n'}),\\
S(n_0,\dots,n_{p+1})=|n_0-n_{p+1}|+(1+\sqrt{n_{p+1}})(1+\sqrt{n'})\sim (1+\sqrt{n_0})^2
\end{gather*} and if $n_0\le Cn_{p+1}$, we have \begin{gather*}\mu(n_0,\dots,n_{p+1})
\le (1+\sqrt{n'})(1+\sqrt{n_{p+1}}),\\
S(n_0,\dots,n_{p+1})\ge c(|n_0-n_{p+1}|)\ge c(n_0+n_{p+1})\sim
(1+\sqrt{n_{p+1}})^2.\end{gather*} In both cases we have \begin{equation} \frac{\mu(n_0,
\dots,n_{p+1})}{S(n_0,\dots,n_{p+1})}\le
C\frac{1+\sqrt{n'}}{1+\sqrt{n_0}+\dots+\sqrt{n_{p+1}}}=C\frac{\max\nolimits_2(\sqrt{n_1},\dots,
\sqrt{n_{p+1}})}{1+\sqrt{n_0}+\dots+\sqrt{n_{p+1}}},\end{equation} where
$\max\nolimits_2(\sqrt{n_1},\dots, \sqrt{n_{p+1}})$ is defined above definition
\ref{def-E3}. Thus theorem \ref{thm-D1} allows us to get (\ref{eqn-E14}) with $\tau=2s$
and some $\nu>0$. (\ref{eqn-E21}) with $\omega=\omega_\ell$ is satisfied by the
definition of $R^{p,1}_\ell$. So $R^{p,1}_\ell\in \widetilde{\mathcal
{R}}^{\nu,2s}_{p+1}(\omega_\ell)$.

\paragraph{}
\textbf{(iii)} The term $I^3_p$.

The treatment of $I^3_p$ is similar to that of $I^2_p$. The only difference is that we
have different support for $B_2$ and $B_3$. So we define
\begin{equation}\label{eqn-F47}R^{p,2}_\ell(u_1,\dots,u_{p+1})
=-\frac12\sum_{n\notin S^\ell_p}\pi_5 \Lambda_m^{2s}\Pi_{n_0}
[(\Pi_{n_1}\Lambda_m^{-1}u_1)\dots(\Pi_{n_{p+1}}\Lambda_m^{-1}u_{p+1})]\end{equation}
with $\pi_5$ given by \begin{equation}\label{eqn-F48}\pi_5=\frac{1}{2^p}\begin{pmatrix}p\\
\ell\end{pmatrix}B_3(n_1,\dots,n_{p+1})\end{equation} and we get
\begin{equation}\label{eqn-F49}I^3_p=\sum_{\ell=0}^pRe\ i\langle R^{p,2}_\ell
(\underset{\ell}{\underbrace{\bar{u},\dots,\bar{u}}},\underset{p+1-\ell}{\underbrace
{u,\dots,u}}),u\rangle.\end{equation} From the support property of $B_3$ we know that
$\Pi_{n_0}R^{p,2}_\ell(\Pi_{n_1}u_1,\dots, \Pi_{n_{p+1}}u_{p+1})$ is supported on domain
$\delta n_{p+1}\le \max\{n_1,\dots,n_p\}=n'\le n_{p+1}$. So on this domain we have
\begin{gather*}\mu(n_0, \dots,n_{p+1})\le (1+\sqrt{n_{p+1}})(1+\sqrt{n'}),\\
S(n_0,\dots,n_{p+1})\sim (1+\sqrt{n_0}+\sqrt{n_{p+1}})^2,\end{gather*} from which we
deduce
\begin{equation}\label{eqn-F50} \frac{\mu(n_0,
\dots,n_{p+1})}{S(n_0,\dots,n_{p+1})}\le
C\frac{1+\sqrt{n'}}{1+\sqrt{n_0}+\dots+\sqrt{n_{p+1}}}.\end{equation} Thus we have by
theorem \ref{eqn-D1}, for any $N\in \N$, there exists $C_N>0$, such that (\ref{eqn-E14})
holds true with $\tau=2s$ and some $\nu>0$. On the other hand, (\ref{eqn-E21}) with
$\omega=\omega_\ell$ is satisfied by  the definition. So $R^{p,2}_\ell\in
\widetilde{\mathcal {R}}^{\nu, 2s}_{p+1}(\omega_\ell)$.
\paragraph{}

Taking $M^p_\ell$ to be the sum of $M^{p,1}_\ell$ and $M^{p,2}_\ell$, and $R^p_\ell$ the
sum of $R^{p,1}_\ell$ and $R^{p,2}_\ell$, we get (\ref{eqn-F7}) with $M^{p}_\ell\in
\widetilde{\mathcal {M}}^{\nu, 2s-a}_{p+1}(\omega_\ell)$ and $R^{p}_\ell\in
\widetilde{\mathcal {R}}^{\nu, 2s}_{p+1}(\omega_\ell)$. This concludes the proof the
lemma.\end{proof}

We have to treat the second term in the right hand side of (\ref{eqn-F6}).
\begin{lem}\label{thm-F3} There are multilinear operators $\widetilde{M}^p_\ell\in
\widetilde{\mathcal {M}}^{\nu, 2s-1}_{p+1}(\widetilde{\omega}_\ell),\
\widetilde{R}^p_\ell \in \widetilde{\mathcal {R}}^{\nu,
2s}_{p+1}(\widetilde{\omega}_\ell)$, $\kappa\le p\le 2\kappa-1$, $ 0\le \ell\le p$ with
$\widetilde{\omega}_\ell$ defined by $\widetilde{\omega}_\ell(j)=-1$, $j=0, \dots, \ell,
p+1$,$\ \widetilde{\omega}_\ell(j)=1, j=\ell+1, \dots, p$, such that
\begin{equation}\label{eqn-F18}\begin{split}Re\ i\langle \Lambda_m^s C(u,\bar{u})\bar{u}, \Lambda_m^s u\rangle
=\sum_{p=\kappa}^{2\kappa-1}\sum_{\ell=0}^{p}&Re\ i\widetilde{M}^p_\ell(\underset{\ell}
{\underbrace{\bar{u}, \dots, \bar{u}}}, \underset{p-\ell}{\underbrace{u, \dots,
u}},\bar{u}), u\rangle\\ &+\sum_{p=\kappa}^{2\kappa-1}\sum_{\ell=0}^{p}Re\ i\langle
\widetilde{R}^p_\ell(\underset{\ell}{\underbrace{\bar{u}, \dots, \bar{u}}},
\underset{p-\ell}{\underbrace{u, \dots, u}}, \bar{u}),
u\rangle.\end{split}\end{equation}\end{lem}
\begin{proof}[Proof of Lemma \ref{thm-F3}:]  Let $\widetilde{\omega}_\ell$ be defined in
the statement of the lemma. We set
\begin{gather}\label{eqn-G2'}\begin{split}
      \widetilde{S}_{p}^\ell=\{(n_0,\dots, n_{p+1})\in\N^{p+2}; \thinspace \mbox{ there exists bijection }
\sigma  \mbox{ from }\\ \{j; 0\le j \le p+1, \widetilde{\omega}_\ell(j)=-1\}  \mbox{ to }
\{ j; 0\le j\le p+1, \widetilde{\omega}_\ell(j)=1\} \\ \mbox{ such that for each } j
\mbox{ in the first set } n_j=n_{\sigma(j)}\}.
   \end{split}\end{gather}
Taking the expression of $C(u,\bar{u})$ defined above (\ref{eqn-F3}) into account, we
compute using notation (\ref{s})
\begin{equation}\label{eqn-F51}\begin{split}&\qquad Re\ i\langle \Lambda_m^{2s}C(u,\bar{u})
\bar{u}, u\rangle\\&=Re\ i\langle
-\frac12\sum_{p=\kappa}^{2\kappa-1}\Lambda_m^{2s}A_p\biggr(\Lambda_m^{-1}(\frac{u+\bar{u}}{2})\biggr)
\Lambda_m^{-1}\bar{u},u\rangle\\&=Re\ i\langle
\sum_{p=\kappa}^{2\kappa-1}\sum_{n\in\N^{p+2}}\sum_{\ell=0}^p
\pi_6\Lambda_m^{2s}\Pi_{n_0}[(\Pi_{n_1}\Lambda_m^{-1}\bar{u})\dots
(\Pi_{n_\ell}\Lambda_m^{-1}\bar{u})\\ &\qquad\qquad\qquad\qquad\qquad\qquad\times
(\Pi_{n_{\ell+1}}\Lambda_m^{-1}u)\dots(\Pi_{n_{p}}\Lambda_m^{-1}u)
(\Pi_{n_{p+1}}\Lambda_m^{-1}\bar{u})],u\rangle
\\&=Re\ i \sum_{p=\kappa}^{2\kappa-1}\sum_{n\in\N^{p+2}}\sum_{\ell=0}^p
\pi_6\int (\Pi_{n_0}\Lambda_m^{2s}\bar{u})(\Pi_{n_1}\Lambda_m^{-1}\bar{u})\dots
(\Pi_{n_\ell}\Lambda_m^{-1}\bar{u})\\
&\qquad\qquad\qquad\qquad\qquad\qquad\times(\Pi_{n_{\ell+1}}\Lambda_m^{-1}u)\dots(\Pi_{n_{p}}\Lambda_m^{-1}u)
(\Pi_{n_{p+1}}\Lambda_m^{-1}\bar{u})dx,\end{split}\end{equation} where $\pi_6$ is given
by
\begin{equation}\label{eqn-F52}\pi_6=-\frac{1}{2^{p+1}}\begin{pmatrix}p\\
\ell\end{pmatrix}B(n_1,\dots,n_{p+1}).\end{equation} With the same reasoning as in the
paragraph above (\ref{eqn-F31}) we may assume $n\notin \widetilde{S}^\ell_p$ in the
computation of (\ref{eqn-F51}). Let $\chi\in C_0^{\infty}(\mathbb{R}), \chi\equiv1$ near
zero, and supp$\chi$ small enough. According to (\ref{eqn-F51}), we define
\begin{equation*}\label{eqn-F19}\begin{split}&\widetilde{M}^p_\ell(u_1, \dots,u_{p+1})=\sum_
{n\notin \widetilde{S}^\ell_p}
\chi\biggr(\frac{|{\lambda_{n_0}^2}-{\lambda_{n_{p+1}}^2}|}{{\lambda_{n_0}^2}+{\lambda_{n_{p+1}}^2}}\biggr)\pi_6
\Lambda_m^{2s}\Pi_{n_0}[(\Pi_{n_1}\Lambda_m^{-1}u_1),\dots,(\Pi_{n_{p+1}}\Lambda_m^{-1}u_{p+1})],
\\& \widetilde{R}^p_\ell(u_1, \dots, u_{p+1})=\sum_{n\notin
\widetilde{S}^\ell_p}\biggr(1-\chi\biggr(\frac
{|\lambda_{n_0}^2-\lambda_{n_{p+1}}^2|}{\lambda_{n_0}^2+\lambda_{n_{p+1}}^2}\biggr)\biggr)\pi_6
\Lambda_m^{2s}\Pi_{n_0}[(\Pi_{n_1}\Lambda_m^{-1}u_1),\dots,(\Pi_{n_{p+1}}\Lambda_m^{-1}u_{p+1})].
\end{split}\end{equation*} It follows that (\ref{eqn-F18}) holds true.

\paragraph{}

Now we are left to check that $\widetilde{M}^{p}_\ell\in \widetilde{\mathcal {M}}^{\nu,
2s-1}_{p+1}(\widetilde{\omega}_\ell)$ and $\widetilde{R}^{p}_\ell\in \widetilde{\mathcal
{R}}^{\nu, 2s}_{p+1}(\widetilde{\omega}_\ell)$.

Because of cut-off function and the support property of function $B$ in the definition of
$\widetilde{M}^p_\ell$ we know that (\ref{eqn-E1}) holds true for $\widetilde{M}^p_\ell$
and we may assume $n_0\sim n_{p+1}$ when estimating $L^2$ norm of
$\Pi_{n_0}\widetilde{M}^p_\ell(\Pi_{n_1}u_1, \dots, \Pi_{n_{p+1}}u_{p+1})$. Since there
is a $\Lambda_m^{-1}$ following each orthogonal projector $\Pi_{n_j}, \quad j=1, \dots,
p+1$, we see that (\ref{eqn-D4}) implies (\ref{eqn-E2}) with $\tau=2s-1$ and some
$\nu>0$. Moreover, (\ref{eqn-E13}) with $\omega=\widetilde{\omega}_\ell$ is satisfied by
the definition of $\widetilde{M}^p_\ell$. So $\widetilde{M}^{p}_\ell\in
\widetilde{\mathcal {M}}^{\nu, 2s-1}_{p+1}(\widetilde{\omega}_\ell)$.

Assume $\Pi_{n_0}[R(\Pi_{n_1}u_1, \dots, \Pi_{n_{p+1}}u_{p+1})]$ does not vanish. Then we
have $|n_0-n_{p+1}|\ge c(n_0+n_{p+1})$ for some small $c>0$ because of the cut-off
function and  $n_{p+1}\ge \max\{n_1,\dots,n_p\}=n'$ because of the support property of
function $B$. Therefore if $n_0\ge n'$, we have
\begin{equation*}\begin{split}\mu(n_0,\dots,n_{p+1})=(1+\sqrt{n'})(1+\min\{\sqrt{n_0},\sqrt{n_{p+1}}\}),\\
S(n_0,\dots,n_{p+1})=|n_0-n_{p+1}|+(1+\sqrt{n'})(1+\min\{\sqrt{n_0},\sqrt{n_{p+1}}\}),\end{split}\end{equation*}
and thus $$\frac{\mu(n_0,\dots,n_{p+1})}{S(n_0,\dots,n_{p+1})}\le
C\frac{1+\sqrt{n'}}{\sqrt{n_0}+\sqrt{n_{p+1}}+1+\sqrt{n'}}\le
C\frac{\max\nolimits_2(\sqrt{n_1},\dots,\sqrt{n_{p+1}})}{1+\sqrt{n_0}+\dots+\sqrt{n_{p+1}}};$$
if $n_0<n'$, we have \begin{equation*}\begin{split} \mu(n_0,\dots,n_{p+1})\le
(1+\sqrt{n'})^2,\qquad
S(n_0,\dots,n_{p+1})=|n'-n_{p+1}|+\mu(n_0,\dots,n_{p+1}),\end{split}\end{equation*} and
thus $$\frac{\mu(n_0,\dots,n_{p+1})}{S(n_0,\dots,n_{p+1})}\le
C\frac{1+\sqrt{n'}}{\sqrt{n'}+\sqrt{n_{p+1}}+1+\sqrt{n'}}\le
C\frac{\max\nolimits_2(\sqrt{n_1},\dots,\sqrt{n_{p+1}})}{1+\sqrt{n_0}+\dots+\sqrt{n_{p+1}}}.$$
Now using theorem \ref{thm-D1} we see that (\ref{eqn-E14}) holds true with $\tau=2s$ and
some $\nu>0$. But (\ref{eqn-E21}) with $\omega=\widetilde{\omega}_\ell$ is satisfied
according to the definition. So $\widetilde{R}^{p}_\ell\in \widetilde{\mathcal {R}}^{\nu,
2s}_{p+1}(\widetilde{\omega}_\ell)$. This concludes the proof of lemma.
\end{proof} Summarizing the above analysis gives an end to the proof of the proposition
\ref{thm-F1}.
\end{proof}

In order to control the energy, let us first turn to some useful estimates in the following subsection.



\subsection{Geometric bounds}   \label{sec-1}
\paragraph{}
This subsection is a modification of section 2.1 in \cite{D}. We give it for the
convenience of the reader. Consider the function on $\mathbb{R}^{p+2}$ depending on the
parameter $m\in(0,+\infty)$, defined for $\ell=0,\dots,p+1$ by
    \begin{equation}\label{eqn-G1}
    F_m^\ell(\xi_0, \dots,\xi_{p+1})=\sum_{j=0}^{\ell}\sqrt{m^2+\xi^2_j}-\sum_{j=\ell+1}^{p+1}\sqrt{m^2+\xi^2_j}.
     \end{equation}

The main result of this subsection is the following theorem:
\begin{thm}\label{thm-G1}
There is a zero measure subset $\mathcal {N}$ of $\R^*_+$ such that for any integers
$0\le\ell\le p+1$, any $m\in\R^*_+-\mathcal{N}$, there are constants $c>
0,N_0\in\mathbb{N}$ such that the lower bound
     \begin{align}\begin{split}\label{eqn-G3}
       |F_m^\ell(\lambda_{n_0},\dots,\lambda_{n_{p+1}})|\ge   c(1+\sqrt{n_0}+\sqrt{n_{p+1}})^{-3-\rho}(1+|\sqrt{n_0}-
       \sqrt{n_{p+1}}|+\sqrt{n'})^{-2N_0}
      \end{split}\end{align}
holds true for any $\rho>0$ and any $(n_0,\dots,n_{p+1})\in \mathbb{N}^{p+2}-S_p^\ell$.
Here $\lambda_n$ are given by (\ref{eqn-B1}), $n'=\max\{n_1, \dots, n_p\}$, and
$S_p^\ell$ is defined in (\ref{eqn-F30}), in which we have set $\omega_\ell(j)=-1,\ j=0,
\dots, \ell,\ \omega(j)=1, \ j=\ell+1, \dots, p+1$.
\end{thm}

The proof of the theorem will rely on some geometric estimates that we shall deduce from
results of \cite{DS1}. Let us show that under the condition of theorem \ref{thm-G1} we
have
\begin{equation}\label{eqn-G4}
 |F_m^\ell(\lambda_{n_0},\dots,\lambda_{n_{p+1}})|\ge   c(1+\sqrt{n_0}+\sqrt{n_{p+1}})^{-3-\rho}(1+
 |\sqrt{n_0}-\sqrt{n_{p+1}}|)
 ^{-N_0}(1+\s{n_1}+\dots+\sqrt{n_p})^{-N_0}.
\end{equation}
Let $I\subset(0,+\infty)$ be some compact interval and define for $0\le\ell\le p+1$ functions
\begin{align}\begin{split}\label{eqn-G5}
f_\ell:[0,1]\times[0,1]^{p+2}\times I &\longrightarrow \mathbb{R}\\
   (z,x_0,\ldots,x_{p+1},y)&\rightarrow f_\ell(z,x_0,\ldots,x_{p+1} ,y)\\
g_\ell:[0,1]\times[0,1]^p\times I &\longrightarrow \mathbb{R}\\
        (z,x_1,\ldots,x_p\ ,y)&\rightarrow g_\ell(z,x_1,\ldots,x_p\ ,y)
\end{split}\end{align}
by \begin{align}\begin{split}\label{eqn-G6}
f_\ell(z,x_0,\ldots,x_{p+1},y)&=\sum_{j=0}^\ell\sqrt{z^2+y^2x_j^2}-\sum_{j=\ell+1}^{p+1}\sqrt{z^2+y^2x_j^2}\\
g_\ell(z,x_1,\ldots,x_p,y)&=z\biggr[\sum_{j=1}^\ell\frac{z}{\sqrt{z^2+y^2x_j^2}}-\sum_{j=\ell+1}^p\frac{z}
{\sqrt{z^2+y^2x_j^2}}
\biggr]\qquad \text{when}\ z>0,\\g_\ell(0,x_1,\ldots,x_p,y)&\equiv 0.
\end{split}\end{align}
Then the graphs of $f_\ell,g_\ell$ are subanalytic subsets of $[0,1]^{p+3}\times I$ and
$[0,1]^{p+1}\times I$ respectively, so that $f_\ell,g_\ell$ are continuous subanalytic
functions (refer to Bierstone-Milman \cite{BM} for an introduction to subanalytic sets
and functions). Let us consider the set $\Gamma$ of points $(z,x)\in[0,1]^{p+3}$(resp.\
$(z,x)\in[0,1]^{p+1})$ such that $y\rightarrow f_\ell(z,x,y)$ (resp. $y\rightarrow
g_\ell(z,x,y)$) vanishes identically. If $(z,x)\in\Gamma$ and $z\ne0$, we have
$$\ell=\frac{p}{2}\qquad \text{and}\qquad \sum_{j\le \ell}x_j^{2\kappa}-\sum_{j\ge\ell+1}x_j^{2\kappa}=0,
\forall \kappa\in\mathbb{N}^*$$
where the sum is taken respectively for $0\le j\le p+1$ in the case of $f_\ell$ and $1\le
j\le p$ for $g_\ell.$ This implies that there is a bijection
$\sigma:\{0,\ldots,\ell\}\rightarrow\{\ell+1,\ldots,p+1\}$ (resp.
$\{1,\ldots,\ell\}\rightarrow\{\ell+1,\ldots,p\}$) such that $x_{\sigma(j)}=x_j$ for any
$j=0,\ldots,\ell$ (resp. $j=1,\ldots,\ell$)---\thinspace see for instance the proof of
lemma 5.6 in \cite{DS1}. When $p$ is even, denote by $S\negmedspace_p$ the set of all
bijections respectively from $\{0,\ldots,\frac p2\}$ to $\{\frac{p}{2}+1,\ldots,p+1\}$
and from $\{1,\ldots,\frac p2\}$ to $\{\frac{p}{2},\ldots,p\}.$ Define for $0\le\ell \le
p+1$
\begin{equation}\begin{split}\label{eqn-G7}
\rho_\ell(z,x)&\equiv z\qquad \mbox{if}\qquad \ell\ne\frac{p}{2},\\
\rho_\ell(z,x)&= z\prod_{\sigma\in S\negmedspace_p}\biggr[\sum_{j\le
p/2}(x_{\sigma(j)}^2-x_j^2)^2 \biggr]\qquad \text{if}\qquad \ell=\frac{p}{2},
\end{split}\end{equation}
where the sum in the above formula is taken for $j\ge0$ (resp. $j\ge 1$) when we study
$f_\ell$ (resp. $g_\ell$). Then the set $\{\rho_\ell=0\}$ contains those points $(z,x)$
such that $y\rightarrow f_\ell(z,x,y)$ (resp. $y\rightarrow g_\ell(z,x,y)$) vanishes
identically. The following proposition is the same as proposition 2.1.2 in \cite{D}.

\begin{pro}\label{thm-G2}(i)There are $\widetilde{N}\in\mathbb{N},
\alpha_0>0,\delta>0,C>0$, such that for any $0\le\ell\le p+1,$ any
$\alpha\in(0,\alpha_0),$ any $(z,x)\in[0,1]^{p+3}$ (resp. $(z,x)\in[0,1]^{p+1}$) with
$\rho_\ell(z,x)\ne0,$ any $N\ge\widetilde{N}$ the sets
\begin{align}\begin{split}\label{eqn-G8}
I^f_\ell(z,x,\alpha)&=\{y\in I;
|f_\ell(z,x,y)|<\alpha\rho_\ell(z,x)^N\}\\
I^g_\ell(z,x,\alpha)&=\{y\in I;|g_\ell(z,x,y)|<\alpha\rho_\ell(z,x)^N\}
\end{split}
\end{align}
have Lebesgue measure bounded from above by $C\alpha^\delta\rho_\ell(z,x)^{N\delta}.$

(ii) For any $N\ge\widetilde{N}$, there is $K\in\mathbb{N}$ such that for any
$\alpha\in(0,\alpha_0)$, any $(z,x)\in[0,1]^{p+1},$ the set $I^g_\ell(z,x,\alpha)$ may be
written as the union of at most $K$ open disjoint subintervals of $I$.
\end{pro}

\paragraph{}
We shall deduce (\ref{eqn-G4}) from several lemmas. Let us first introduce some
notations. When $p$ is odd or $p$ is even and $\ell\ne\frac{p}{2}$, we set
$\mathbb{N}'^p_\ell=\emptyset$. When $p$ is even and $\ell=\frac p2,$ we define
\begin{align}\begin{split}\label{eqn-G9}
\mathbb{N}'^p_\ell=\{\tilde{n}=(n_1,\dots,n_p)\in\mathbb{N}^p;\ &\text{there is a bijection}\ \
\\\sigma:\{1,\dots,\ell\}\rightarrow\{\ell+1,\dots,p\}\ &\text{such
that}\ n_{\sigma(j)}=n_j,j=1,\dots,\ell\}.
\end{split}\end{align}
We set also
\begin{equation}\label{eqn-G10}\mathbb{N}^{p+2}_\ell
=\{(n_0,\dots,n_{p+1})\in\mathbb{N}^{p+2}; \ \tilde{n}\in\mathbb{N}'^p_\ell\ \mbox{and}\ n_0=n_{p+1}\}.
\end{equation} Of course,
$\mathbb{N}^{p+2}_\ell=\emptyset$ if $p$ is odd or $p$ is even and $\ell\ne\frac p2$.\paragraph{}

We remark first that it is enough  to prove (\ref{eqn-G4}) for those $(n_1,\dots,n_p)$
which do not belong to $\mathbb{N}'^p_\ell$: actually if $p$ is even, $\ell=\frac p2$ and
$(n_1,\dots,n_p)\in\mathbb{N}'^p_\ell$, we have
$|F^\ell_m(\lambda_{n_0},\dots,\lambda_{n_{p+1}})|=|\sqrt{m^2+\lambda_{n_0}^2}-\sqrt{m^2+\lambda_{n_{p+1}}^2}|$
which is bounded from below, when $m$ stays in some compact interval, by
$$\frac{2|n_0-n_{p+1}|}{\sqrt{m^2+\lambda_{n_0}^2}+\sqrt{m^2+\lambda_{n_{p+1}}^2}}\ge
\frac{c}{1+\lambda_{n_0}+\lambda_{n_{p+1}}}$$ since from
$(n_0,\dots,n_{p+1})\in\mathbb{N}^{p+2}-S_p^\ell$, we have $n_0\ne n_{p+1}$. Consequently
(\ref{eqn-G4}) holds true trivially. From now on, we shall always consider $p-$tuple
$\tilde{n}$ which do not belong to $\mathbb{N}'^p_\ell$.

Let us define for $\ell=1,\dots,p$ another function on $\mathbb{R}^p$ given by
\begin{equation}\label{eqn-G11}
G^\ell_m(\xi_1,\dots,\xi_p)=\sum_{j=1}^\ell\sqrt{m^2+\xi_j^2}-\sum_{j=\ell+1}^p\sqrt{m^2+\xi_j^2}.\end{equation}
Let $J\subset(0,+\infty)$ be a given compact interval. For $\alpha>0, N_0\in\mathbb{N},
0\le\ell\le p+1$, $n=(n_0,\dots,n_{p+1})\in\mathbb{N}^{p+2}$ define
\begin{align}\begin{split}\label{eqn-G12}
E_J^\ell(n,\alpha,N_0)=\{m\in
J;|F_m^\ell(\lambda_{n_0},\dots,\lambda_{n_{p+1}})|<\alpha(1+\lambda_{n_0}+\lambda_{n_{p+1}})^{-3-\rho}\\
\times(1+|\lambda_{n_0}-\lambda_{n_{p+1}}|)^{-N_0}(1+\lambda_{n_1}+\dots+\lambda_{n_p})^{-N_0}\}.
\end{split}\end{align}
We set also for $\beta>0,N_1\in\mathbb{N}^*, \tilde{n}=(n_1,\dots,n_p)\in\mathbb{N}^p-\mathbb{N}'^p_\ell$
\begin{equation}\label{eqn-G13}
E'^\ell_J(\tilde{n},\beta,N_1)=\{m\in J;\biggr|\frac{\partial G_m^\ell}{\partial
m}(\lambda_{n_1},\dots,\lambda_{n_p})\biggr|<\beta(1+\lambda_{n_1}+\dots+\lambda_{n_p})^{-N_1}\}.
\end{equation}
We define for $\gamma>\beta$ a subset of $\mathbb{N}^{p+2}$ by
\begin{align}\begin{split}\label{eqn-G14}
S(\beta,\gamma,N_1)=\{(n_0,\dots,n_{p+1})\in\mathbb{N}^{p+2}&-\mathbb{N}^{p+2}_\ell:
\lambda_{n_0}<\frac{\gamma}{3\beta}(1+\lambda_{n_1}+\dots+\lambda_{n_p})^{N_1}\\
&\mbox{or}\qquad \lambda_{n_{p+1}}<\frac{\gamma}{3\beta}(1+\lambda_{n_1}+\dots+\lambda_{n_p})^{N_1}\}.
\end{split}\end{align}

\begin{lem}\label{thm-G3} Let $\widetilde{N}, \delta,\alpha_0$ be the constants defined in the
statement of proposition \ref{thm-G2}. There are constants $C_1>0, M\in\mathbb{N}^*$ such
that for any $\beta\in(0,\alpha_0)$, any $N_1\in\mathbb{N}$ with $N_1>M\widetilde{N}$ and
$N_1>\frac{2pM}{\delta}$, one has
\begin{equation}\label{eqn-G15}
meas\biggr[\bigcup_{\tilde{n}\in\mathbb{N}^p-\mathbb{N}'^p_\ell}E'^\ell_J(\tilde{n},\beta,N_1)\biggr]\le
C_1\beta^\delta.\end{equation}
\end{lem}
\begin{proof}
Set $y=\frac1m$ and $$z=(1+\sum_{j=1}^p\lambda_{n_j})^{-1},\ x_j=\lambda_{n_j}z\ ,\
j=1,\dots,p.$$ Denote by $X$ the set of points $(z,x)\in[0,1]^{p+1}$ of the preceding
form for $(n_1, \dots,n_p)$ describing $\mathbb{N}^p$. When $p$ is even and $\ell=p/2$,
let $X'^p_\ell$ be the imagine of $\mathbb{N}'^p_\ell$ defined by (\ref{eqn-G9}) under
the map $\tilde{n}\rightarrow(z,x)$. Using definition (\ref{eqn-G7}), we see that there
are constants $M>0, C>0$, depending only on $p$, such that for $0\le\ell\le
p+1$\begin{equation}\label{eqn-G16} \forall(z,x)\in X-X'^p_\ell\ ,\
z^M\le\rho_\ell(z,x)\le Cz\end{equation} since , when $\ell=\frac{p}{2}$ and $(n_1,\dots,
n_p)\notin \mathbb{N}'^p_\ell$, $\sum_{j=1}^{\frac{p}{2}}
(\lambda^2_{n_{\sigma(j)}}-\lambda_{n_j}^2)^2\ge 1$, by the definition of
$\lambda_{n_j}$. Remark that  with the above notations
$$\frac{\partial G^\ell_m}{\partial m}(\lambda_{n_1},\dots,\lambda_{n_p})=\sum_{j=1}^\ell
\frac{m}{\sqrt{m^2+\lambda_{n_j}^2}}-\sum_{j=\ell+1}^p
\frac{m}{\sqrt{m^2+\lambda_{n_j}^2}}=\frac1zg_\ell(z,x_1,\dots,x_p,y).$$Then if
$I=\{m^{-1}; m\in J\},$ we see that $m\in E'^\ell_J(\tilde{n},\beta,N_1)$ for
$n\notin\mathbb{N}'^p_\ell$ if and only if $y=\frac1m$ satisfies
\begin{equation}\label{eqn-G17}
|g_\ell(z,x_1,\dots,x_p,y)|<\beta z^{N_1+1}\le\beta\rho_\ell(z,x)^{\frac1M(N_1+1)}
\end{equation} using (\ref{eqn-G16}). Applying proposition \ref{thm-G2} (i),
we see that for any fixed value of $(z,x)\in X-X'^p_\ell$, the measure of those $y$ such
that (\ref{eqn-G17}) holds true is bounded from above by
$$C\beta^\delta\rho_\ell(z,x)^{\frac{N_1+1}{M}\delta}\le C\beta^\delta z^{\frac{N_1+1}{M}\delta}$$
if we assume $N_1\ge M\widetilde{N}$ and $\beta\in(0,\alpha_0)$. Consequently, we get
with a constant $C'$ depending only on $J$,
\begin{align*}meas(E'^\ell_J(n',\beta,N_1))\le
C'\beta^\delta(1+\lambda_{n_1}+\dots+\lambda_{n_p})^{-\frac{N_1+1}{M}\delta}\\\le
C'\beta^\delta(1+n_1+\dots+n_p)^{-\frac{N_1+1}{2M}\delta}.\end{align*} Inequality
(\ref{eqn-G15}) follows from this estimate and the assumption on $N_1.$
\end{proof}

\begin{lem}\label{thm-G4}
Let $\widetilde{N}, \delta,\alpha_0$ be the constants defined in the statement of
proposition \ref{thm-G2}. There are constants $M\in\mathbb{N}^*, \theta>1, C_2>0$ such
that for any $N_0,N_1\in\mathbb{N}^*$ satisfying $N_0>\widetilde{N}MN_1$ and
$N_0\delta>2(p+2)MN_1$, any $0<\beta<\gamma$ with $\frac{\gamma}{\beta}>\theta$, any
$\alpha>0$ satisfying $\alpha(\frac{\beta}{2\gamma})^{-\frac{N_0}{N_1}}<\alpha_0$, one
has
\begin{equation}\label{eqn-G18}
meas\biggr[\bigcup_{n\in S(\beta,\gamma,N_1)}E^\ell_J(n,\alpha,N_0)\biggr]\le
C_2\alpha^\delta\biggr(\frac{\beta}{2\gamma}\biggr)^{-\frac{N_0}{N_1}\delta}.\end{equation}
\end{lem}
\begin{proof} We first remark that if $\lambda_{n_0}+\lambda_{n_{p+1}}>\frac\gamma\beta
(1+\lambda_{n_1}+\dots+\lambda_{n_p})^{N_1}$ and  $n\in S(\beta,\gamma,N_1)$, then either
$$ \lambda_{n_0}\ge\frac{2\gamma}{3\beta}(1+\lambda_{n_1}+\dots+\lambda_{n_p})^{N_1}\qquad
 \mbox{or}\qquad
\lambda_{n_{p+1}}\ge\frac{2\gamma}{3\beta} (1+\lambda_{n_1}+\dots+\lambda_{n_p})^{N_1},$$
which implies that
$$|F_m^\ell(\lambda_{n_0},\dots,\lambda_{n_{p+1}})|\ge c\frac\gamma\beta(1+\lambda_{n_1}
+\dots+\lambda_{n_p})^{N_1}$$ for some constant $c>0$ depending only on $p$ and $J$, if
$\frac\gamma\beta>\theta$ large enough. Consequently, if $\alpha<\alpha_0$ small enough
relatively to $c$, we see that we have in this case $E^\ell_J(n,\alpha,N_0)=\emptyset$
when $n\in S(\beta,\gamma,N_1)$. We may therefore consider only indices $n$ such that $$
n\in S(\beta,\gamma,N_1)\qquad \mbox{and}\qquad
\lambda_{n_0}+\lambda_{n_{p+1}}\le\frac\gamma\beta(1+\lambda_{n_1}+\dots+\lambda_{n_p})^{N_1}.$$
Consequently, for $m\in E^\ell_J(n,\alpha,N_0)$ and $n\in S(\beta,\gamma,N_1)$, we have
\begin{align}\begin{split}\label{eqn-G19}
|F_m^\ell(\lambda_{n_0},\dots,\lambda_{n_{p+1}})|&\le\alpha(1+\lambda_{n_1}+\dots+\lambda_{n_p})^{-N_0}\\
&\le\alpha\biggr(\frac{\beta}{2\gamma}\biggr)^{-\frac{N_0}{N_1}}(1+\lambda_{n_0}+\dots+
\lambda_{n_{p+1}})^{-\frac{N_0}{N_1}}.\end{split}\end{align} Define for
$n\in\mathbb{N}^{p+2}$
\begin{equation}\label{eqn-G20}z=(1+\sum_{j=0}^{p+1}\lambda_{n_j})^{-1},\qquad
x_j=\lambda_{n_j}z,\  j=0, \dots, p+1.\end{equation}

Denote by $X\subset[0,1]^{p+3}$ the set of points $(z,x)$ of the preceding form, and let
$X^p_\ell$ be the imagine of the set $\mathbb{N}_\ell^{p+2}$ defined by (\ref{eqn-G10})
under the map $n\rightarrow(z,x)$. By (\ref{eqn-G7}) we have again
$$\forall(z,x)\in X-X^p_\ell\ , \ z^M\le\rho_\ell(z,x)\le Cz$$
for some large enough $M$, depending only on $p$. Moreover$$F^\ell_m(\lambda_{n_0},\dots,
\lambda_{n_{p+1}})=\frac mzf_\ell(z,x_0,\dots,x_{p+1},y)$$ and (\ref{eqn-G19}) implies
that if $n\in S(\beta,\gamma,N_1)$ and $m\in E^\ell_J(n,\alpha,N_0)$, then $y$ satisfies
\begin{align}\begin{split}\label{eqn-G21}
|f_\ell(z,x_0,\dots,x_{p+1}, y)|&\le
C\alpha\biggr(\frac{\beta}{2\gamma}\biggr)^{-\frac{N_0}{N_1}}z^{1+\frac{N_0}{N_1}}\\
&\le C\alpha\biggr(\frac{\beta}{2\gamma}\biggr)^{-\frac{N_0}{N_1}}\rho_\ell(z,x)^
{\frac1M(1+\frac{N_0}{N_1})}\end{split}\end{align} We assume that $\alpha,N_0,N_1$
satisfy the conditions of the statement of the lemma. Then by (i) of proposition
\ref{thm-G2} we get that the measure of those $y\in J$ satisfying (\ref{eqn-G21}) is
bounded from above by
$$C\biggr[\alpha\biggr(\frac{\beta}{2\gamma}\biggr)^{-\frac{N_0}{N_1}}\biggr]^\delta z^
{\frac\delta M(1+\frac{N_0}{N_1})}$$
for some constant $C$, independent of $N_0,N_1,\alpha,\beta,\gamma$. Consequently the
measure of $E^\ell_J(n,\alpha,N_0)$ is bounded from above when $n\in S(\beta,\gamma,N_1)$
by
\begin{multline*}C\biggr[\alpha\biggr(\frac{\beta}{2\gamma}\biggr)^{-\frac{N_0}{N_1}}\biggr]
^\delta\biggr(1+\lambda_{n_0}+\dots+\lambda_{n_{p+1}}\biggr) ^{-\frac\delta
M(1+\frac{N_0}{N_1})}\\ \le
C'\biggr[\alpha\biggr(\frac{\beta}{2\gamma}\biggr)^{-\frac{N_0}{N_1}}\biggr]^\delta\biggr(1+n_0
+\dots+n_{p+1}\biggr)^{-\frac{\delta}{2M}(1+\frac{N_0}{N_1})}\end{multline*} for another
constant $C'$ depending on $J$. The conclusion of the lemma follows by summation, using
that $\frac{\delta}{M}(1+\frac{N_0}{N_1})>2(p+2)$.\end{proof}
 \paragraph{}

\begin{proof}[Proof of theorem \ref{thm-G1}]: We fix $N_0,N_1$ satisfying the conditions
stated in lemmas \ref{thm-G3} and \ref{thm-G4}, and such that $N_0>2p+N_1$. We write when
$n\notin S(\beta,\gamma,N_1),\ 0\le\ell\le p+1$,$$
E^\ell_J(n,\alpha,N_0)\subset[E^\ell_J(n,\alpha,N_0)\cap
E'^\ell_J(\tilde{n},\beta,N_1)]\cup[E^\ell_J(n,\alpha,N_0)\cap(E'^\ell_J(\tilde{n},\beta,N_1)^c]$$
and estimate, using that we reduced ourselves to those $\tilde{n}\notin
\mathbb{N}'^p_\ell$
\begin{align}\begin{split}\label{eqn-G22}
meas\biggr[\bigcup_{n;\ \tilde{n}\notin\mathbb{N}'^p_\ell}E^\ell_J(n,\alpha,N_0)\biggr]
&\le meas\biggr[\bigcup_{n\in
S(\beta,\gamma,N_1)}E^\ell_J(n,\alpha,N_0)\biggr]+meas\biggr[\bigcup_{\tilde{n}\notin\mathbb{N}'^p_\ell}
E'^\ell_J(\tilde{n},\beta,N_1)\biggr]\\
&+meas\biggr[\bigcup_{n\in
S(\beta,\gamma,N_0)^c-\mathbb{N}^{p+2}_\ell}E^\ell_J(n,\alpha,N_0)\cap
E'^\ell_J(\tilde{n},\beta,N_1)^c\biggr].\end{split}\end{align} Let us bound the measure
of $E^\ell_J(n,\alpha,N_0)\cap E'^\ell_J(\tilde{n},\beta,N_1)^c$ for $n\in
S(\beta,\gamma,N_0)^c-\mathbb{N}^{p+2}_\ell$. If $m$ belongs to that set, the inequality
in (\ref{eqn-G12}) holds true. Remark that we may assume $\ell\le p:$ if $\ell=p+1,$
$|F^\ell_m(\lambda_{n_0},\dots,\lambda_{n_{p+1}})|\ge
c(1+\lambda_{n_0}+\lambda_{n_{p+1}})$ for some $c>0$, which is not compatible with
(\ref{eqn-G12}) for $\alpha<\alpha_0$ small enough. Let us write (\ref{eqn-G12}) as
\begin{equation}\begin{split}\label{eqn-G23}
|\lambda_{n_0}-\lambda_{n_{p+1}}+\widetilde{G}_m(\lambda_{n_0},\dots,\lambda_{n_{p+1}})|
<\alpha(1+\lambda_{n_0}+\lambda_{n_{p+1}})^{-3-\rho}
\\ \times(1+|\lambda_{n_0}-\lambda_{n_{p+1}}|)^{-N_0}(1+\lambda_{n_1}+\dots+\lambda_{n_p})
^{-N_0}\end{split}\end{equation}
with, using notation (\ref{eqn-G11})
\begin{align}\begin{split}\label{eqn-24}
\widetilde{G}_m(\lambda_{n_0},\dots,\lambda_{n_{p+1}})&=G_m(\lambda_{n_1},\dots,\lambda_{n_{p}})
+R_m(\lambda_{n_0},\lambda_{n_{p+1}})\\
R_m(\lambda_{n_0},\lambda_{n_{p+1}})&=(\sqrt{m^2+\lambda_{n_0}^2}-\lambda_{n_0})-(\sqrt{m^2
+\lambda_{n_{p+1}}^2}-\lambda_{n_{p+1}}).
\end{split}\end{align}
Since $n\in S(\beta,\gamma,N_1)^c$, we have by (\ref{eqn-G14})
\begin{equation}\label{eqn-G25}
\lambda_{n_0}\ge\frac{\gamma}{3\beta}(1+\lambda_{n_1}+\dots+\lambda_{n_p})^{N_1}\ ,\
\lambda_{n_{p+1}}\ge\frac{\gamma}{3\beta}(1+\lambda_{n_1}+\dots+\lambda_{n_p})^{N_1}.
\end{equation}
Consequently there is a constant $C>0$, depending only on $J$, such that $$|\frac{\partial R_m}{\partial
m}(\lambda_{n_0},\lambda_{n_{p+1}})|\le C\frac{\beta}{\gamma}(1+\lambda_{n_1}+\dots+\lambda_{n_p})
^{-N_1}.$$ If $\gamma$ is large
enough and $m\in E'^\ell_J(\tilde{n},\beta,N_1)^c$, we deduce from (\ref{eqn-G13}) that
\begin{equation}\label{eqn-G26}
|\frac{\partial\widetilde{G_m}}{\partial
m}(\lambda_{n_0},\dots,\lambda_{n_{p+1}})|\ge\frac\beta2(1+\lambda_{n_1}+\dots+\lambda_{n_p})^{-N_1}.
\end{equation}
By (ii) of proposition \ref{thm-G2}, we know that there is $K\in\mathbb{N}$, independent
of $\alpha, \beta,\gamma$ such that the set $J-E'^\ell_J(\tilde{n},\beta,N_1)$ is the
union of at most $K$ disjoint intervals $J_j(\tilde{n},\beta,N_1)\ ,\ 1\le j\le K$.
Consequently, we have
\begin{equation}\label{eqn-G27}
E^\ell_J(n,\alpha,N_0)\cap (E'^\ell_J(\tilde{n},\beta,N_1))^c\subset\bigcup_{j=1}^K\{m\in J_j
(\tilde{n},\beta,N_1);\
(\ref{eqn-G23})\ \mbox{holds true}\},
\end{equation}
and on each interval $J_j(n',\beta,N_1)$, (\ref{eqn-G26}) holds true. We may on each such
interval perform in the characteristic function of (\ref{eqn-G23}) the change of variable
of integration given by
$m\rightarrow\widetilde{G}_m(\lambda_{n_0},\dots,\lambda_{n_{p+1}})$. Because of
(\ref{eqn-G26}) this allows us to estimate the measure of (\ref{eqn-G27}) by
\begin{align*}K\frac2\beta\alpha(1+\lambda_{n_0}+\lambda_{n_{p+1}})^{-3-\rho}(1+|\lambda_{n_0}-\lambda_{n_{p+1}}|)^{-N_0}
(1+\lambda_{n_1}+\dots+\lambda_{n_p})^{-N_0+N_1}\\ \le
CK\frac2\beta\alpha(1+n_0+n_{p+1})^{-\frac12(3+\rho)}(1+|\sqrt{n_0}-\sqrt{n_{p+1}}|)^{-N_0}
(1+n_1+\dots+n_p)^{-\frac12(N_0-N_1)}\end{align*} Summing in $n_0,\dots,n_{p+1}$, we see
that since $N_0>2p+N_1$, the last term in (\ref{eqn-G22}) is bounded from above by
$C_3\frac\alpha\beta$ with $C_3$ independent of $\alpha,\beta,\gamma$. By lemmas
\ref{thm-G3} and \ref{thm-G4}, we may thus bound (\ref{eqn-G22}) by $$
C_2\alpha^\delta\biggr(\frac{\beta}{2\gamma}\biggr)^{-\frac{N_0}{N_1}\delta}+C_1\beta^\delta+C_3\frac\alpha\beta$$
if $\alpha,\beta$ are small enough, $\gamma$ is large enough and
$\alpha(\frac{\beta}{\gamma})^{-\frac{N_0}{N_1}}$ is small enough. If we take
$\beta=\alpha^\sigma,\ \gamma=\alpha^{-\sigma}$ with $\sigma>0$ small enough, and
$\alpha\ll 1$, we finally get for some $\delta'>0$,
$$meas\biggr[\bigcup_{n;\
 \tilde{n}\notin\mathbb{N}'^p_\ell}E^\ell_J(n,\alpha,N_0)\biggr]\le C\alpha^{\delta'}\rightarrow0
 \ \mbox{if}\ \alpha\rightarrow
0^+.$$ This implies that in this case the  set of those $m\in J$ for which (\ref{eqn-G4})
does not hold true for any $c>0$ is of
 zero measure. This concludes the proof. \end{proof}
 \paragraph{}
We will need a consequence of theorem \ref{thm-G1}:.
\begin{pro}\label{thm-G6}
There is a zero measure subset $\mathcal {N}$ of $\R^*_+$ such that for any integers
$0\le\ell\le p+1$, any $m\in\R^*_+-\mathcal{N}$, there are constants $c>
0,N_0\in\mathbb{N}$ such that the lower bound
     \begin{align}\begin{split}\label{eqn-G32}
       |F_m^\ell(\lambda_{n_0},\dots,\lambda_{n_{p+1}})|\ge   c(1+\sqrt{n_0}+\sqrt{n_{p+1}})
       ^{-3-\rho}(1+\sqrt{n'})^{-2N_0}
       \frac{\mu(n_0,\dots, n_{p+1})^{2N_0}}{S(n_0,\dots, n_{p+1})^{2N_0}}
      \end{split}\end{align}
holds true for any $\rho>0$ and any $(n_0,\dots,n_{p+1})\in \mathbb{N}^{p+2}-S_p^\ell$
with $n_0\sim n_{p+1}$ and $n_{p+1}\ge n'$. Here $\lambda_n$, $n'$, $S^\ell_p$ are the
same as those in theorem \ref{thm-G1}.
\end{pro}
\begin{proof}By theorem \ref{thm-G1} we know (\ref{eqn-G3}) holds true under the conditions of the proposition.
Since we assume $n_0\sim n_{p+1}$ and $n_{p+1}\ge n'$, we have by (\ref{eqn-D2}) and
(\ref{eqn-D3}) \begin{align}\nonumber\mu(n_0, \dots, n_{p+1})\sim
(1+\sqrt{n_{p+1}})(1+\sqrt{n'}),\\  S(n_0,\dots, n_{p+1})\sim
|n_0-n_{p+1}|+(1+\sqrt{n_{p+1}})(1+\sqrt{n'})\\ \nonumber\sim
(1+\sqrt{n_{p+1}})(1+|\sqrt{n_0}-\sqrt{n_{p+1}}|+\sqrt{n'}).\end{align} Therefore we
deduce from (\ref{eqn-G3}) \begin{equation*}\begin{split}
|F_m^\ell(\lambda_{n_0},\dots,\lambda_{n_{p+1}})| &\ge
c(1+\sqrt{n_0}+\sqrt{n_{p+1}})^{-3-\rho} \frac{(1+\sqrt{n_{p+1}})^{2N_0}}{S(n_0,\dots,
n_{p+1})^{2N_0}}\\& \ge c(1+\sqrt{n_0}+\sqrt{n_{p+1}})^{-3-\rho}(1+\sqrt{n'})^{-2N_0}
\frac{\mu(n_0,\dots, n_{p+1})^{2N_0}}{S(n_0,\dots,
n_{p+1})^{2N_0}}.\end{split}\end{equation*} This concludes the proof of the proposition.
\end{proof}

In the following subsection, we shall also use a simpler version of theorem \ref{thm-G1}.
Let us introduce some notations. For $m\in \R^*_+$, $\xi_j\in\mathbb{R}$, $j=0, \dots,
p+1$, $e=(e_0, \dots, e_{p+1})\in \{-1, 1\}^{p+2}$, define
\begin{equation}\label{eqn-G29} \widetilde{F}_m^{(e)}(\xi_0, \dots, \xi_{p+1})=\sum
_{j=0}^{p+1}e_j\sqrt{m^2+\xi^2_j}.\end{equation} When $p$ is even and
$\sharp\{j;e_j=1\}=\frac{p}{2}+1$, denote by $N^{(e)}$ the set of all $(n_0, \dots,
n_{p+1})\in \mathbb{N}^{p+2}$ such that there is a bijection $\sigma$ from $\{j; 0\le
j\le p+1, e_j=1\}$ to $\{j; 0\le j\le p+1, e_j=-1\}$ so that for any $j$ in the first set
$n_j=n_{\sigma(j)}$. In the other cases, set $N^{(e)}=\emptyset$.
\begin{pro}\label{thm-G5} There is a zero measure subset $\mathcal {N}$ of $\R^*_+$ and for any
$m\in\R^*_+-\mathcal{N}$, there are constants $c> 0,N_0\in\mathbb{N}$ such that for any
$(n_0,\dots,n_{p+1})\in \mathbb{N}^{p+2}-N^{(e)}$ one has
\begin{equation}\label{eqn-G30}
       |\widetilde{F}_m^{(e)}(\lambda_{n_0},\dots,\lambda_{n_{p+1}})|\ge   c(1+\sqrt{n_0}+\dots
       +\sqrt{n_{p+1}})^{-N_0}.
\end{equation} Moreover, if $e_0e_{p+1}=1$, one has the inequality
\begin{equation}\label{eqn-G31}|\widetilde{F}_m^{(e)}(\lambda_{n_0},\dots,\lambda_{n_{p+1}})|
\ge c(1+\sqrt{n_0}+\sqrt{n_{p+1}})
(1+\sqrt{n_1}+\dots+\sqrt{n_p})^{-N_0}.\end{equation}\end{pro} \begin{proof} With the
reasoning as in the proof of proposition 2.1.5 in \cite{D}, we get just by replacing
$(n_0, \dots, n_{p+1})$ with $(\lambda_0, \dots, \lambda_{p+1})$
\begin{equation*}|\widetilde{F}_m^{(e)}(\lambda_{n_0},\dots,\lambda_{n_{p+1}})|\ge
c(1+\lambda_{n_0}+\dots+\lambda_{n_{p+1}})^{-N_0}\end{equation*} and
\begin{equation*} |\widetilde{F}_m^{(e)}(\lambda_{n_0},\dots,\lambda_{n_{p+1}})|\ge
c(1+\lambda_{n_0}+\lambda_{n_{p+1}})(1+\lambda_{n_1}+\dots+\lambda_{n_{p}})^{-N_0}\end{equation*}
when $e_0e_{p+1}=1$. This concludes the proof of the proposition by noting
(\ref{eqn-B1}).
\end{proof}



\subsection{Energy control and proof of main theorem}
We shall use the results of subsection \ref{sec-1} to control the energy. When $M(u_1,
\dots, u_{p+1})$ is a $p+1$-linear form, let us define for $0\le \ell\le p+1$,
\begin{align}\label{eqn-H1}   L^{-}_\ell(M)(u_1, \dots, _{p+1})=-
\Lambda_mM(u_1, \dots, u_{p+1})\\ - \nonumber\sum_{j=1}^\ell M(u_1, \dots, \Lambda_mu_j,
\dots, u_{p+1})+\sum_{j=\ell+1}^{p+1}M(u_1, \dots, \Lambda_mu_{j}, \dots, u_{p+1})
\end{align} and
\begin{align}\label{eqn-H1'}L^{+}_\ell(M)(u_1, \dots, _{p+1})=-
\Lambda_mM(u_1, \dots, u_{p+1})-\sum_{j=1}^\ell M(u_1, \dots, \Lambda_mu_j, \dots,
u_{p+1})\\ \nonumber+\sum_{j=\ell+1}^{p}M(u_1, \dots, \Lambda_mu_{j}, \dots,
u_{p+1})-M(u_1, \dots, u_p, \Lambda_mu_{p+1}).\end{align} We shall need the following
lemma:

\begin{lem}\label{thm-H1} Let $\mathcal {N}$ be the zero measure subset of $\R^*_+$ defined
by taking the union of the zero measure subsets defined in proposition \ref{thm-G6} and
proposition \ref{thm-G5}, and fix $m\in\R^*_+-\mathcal{N}$. Let $\omega_\ell,
\widetilde{\omega}_\ell$ be defined in the statement of proposition \ref{thm-F1}. There
is a $\bar{\nu}\in \mathbb{N}$ such that the following statements hold true for any large
enough integer $s$, any integer $p$ with
  $\kappa\le p\le 2\kappa-1$, any integer $\ell$ with $0\le \ell \le p$, any $\rho>0$:
      \begin{itemize}
\item Let $\theta \in (0,1)$, $M_\ell^p\in \widetilde{\mathcal{M}}_{p+1}^{\nu,2s-a}(\omega_\ell)$ with $a=2$ if $d\ge 2$
and $a=\frac{13}{6}-\varsigma$ for any $\varsigma\in (0,1)$ if $d=1$  and
$\widetilde{M}_\ell^p\in \widetilde{\mathcal{M}}_{p+1}^
{\nu,2s-1}(\widetilde{\omega}_\ell)$. Define \begin{equation}\label{eqn-H2}
M_\ell^{p,\epsilon}(u_1, \dots, u_{p+1})=\sum_{n_0}\sum_{n_{p+1}}
\textbf{1}_{\{\sqrt{n_0}+\sqrt{n_{p+1}}<\epsilon^{-\theta\kappa}\}}\Pi_{n_0}M^p_\ell(u_1,
\dots, u_p, \Pi_{n_{p+1}}u_{p+1}).\end{equation} Then there are $\underline {M}
_\ell^{p,\epsilon}\in\widetilde {\mathcal{M}}_{p+1}^{\nu+\bar{\nu},2s-1}(\omega_\ell)$
and $\underline {M}_\ell^p\in\widetilde{\mathcal{M}}_{p+1}^{\nu,2s-2}
(\widetilde{\omega}_\ell)$ satisfying
\begin{equation}\label{eqn-H3}\begin{split}L_\ell^-(\underline {M}_\ell^{p,\epsilon})(u_1,\dots, u_{p+1})
=M_{\ell}^{p,\epsilon}(u_1,\dots,u_{p+1}), \\ L_\ell^+(\underline
{M}_\ell^{p})(u_1,\dots, u_{p+1})=\widetilde{M}_\ell^p(u_1, \dots,
u_{p+1})\end{split}\end{equation}   with the estimate for all $N\ge \bar{\nu}$,
\begin{equation}\label{eqn-H4}\begin{split}||\underline {M}_{\ell}^{p,\epsilon}||_{\mathcal{M}_{p+1, N}^
{\nu+\bar{\nu},2s-1}}\le C\epsilon^ {-(4-a+\rho)\theta\kappa}
||M_\ell^p||_{\mathcal{M}_{p+1,N}^{\nu, 2s-a}}, \\ ||\underline
{M}_{\ell}^{p}||_{\mathcal{M}_{p+1, N}^{\nu+\bar{\nu},2s-2}}\le
C||\widetilde{M}_\ell^p||_{\mathcal{M}_{p+1,N}^ {\nu,2s-1}},
\end{split}\end{equation} where $||\cdot||_{\mathcal{M}_{p+1, N}^{\nu,\tau}}$ is defined in the statement
of definition \ref{def-E1}.
\item Let $R_\ell^p\in \widetilde{R}_{p+1}^{\nu,2s}(\omega_\ell), \widetilde{R}_\ell^p\in\widetilde{R}_{p+1}^
{\nu,2s}(\widetilde{\omega}_\ell)$.
 Then there are
      $\underline {R}_\ell^p\in\widetilde{\mathcal{R}}_{p+1}^{\nu+\bar{\nu},2s}(\omega_\ell)$ and
      $\underline {R}'^p_\ell\in\widetilde
      {\mathcal{R}}_{p+1}^
      {\nu+\bar{\nu},2s}(\widetilde{\omega}_\ell)$ such that \begin{equation}\label{eqn-H5}\begin{split}
      L_\ell^-(\underline {R}_\ell^p)(u_1,\dots,u_{p+1})=R_\ell^p(u_1,\dots,u_{p+1}), \\ L_\ell^+
      (\underline {R}'^p_\ell)(u_1,\dots,
      u_{p+1})=\widetilde{R}_\ell^p(u_1,\dots,u_{p+1}). \end{split}\end{equation}
      \end{itemize}
  \end{lem}
\begin{proof}(i) We substitute in (\ref{eqn-H3}) $\Pi_{n_j}u_j$ to $u_j, j=1, \dots, p+1$, and
compose on the left with $\Pi_{n_0}$. According to (\ref{eqn-H1}), equalities in
(\ref{eqn-H3}) may be written
\begin{align}\label{eqn-H6}
-F_m^\ell(\lambda_{n_0}, \dots, \lambda_{n_{p+1}})\Pi_{n_0}\underline
{M}_\ell^{p,\epsilon}(\Pi_{n_1}u_1,\dots, \Pi_{n_{p+1}}u_{p+1})=\Pi_{n_0}M_\ell^{p,
\epsilon}(\Pi_{n_1}u_1,\dots, \Pi_{n_{p+1}}u_{p+1}),
\\ \label{eqn-H6'} \widetilde{F}_m^{(e)}(\lambda_{n_0}, \dots, \lambda_{n_{p+1}})\Pi_{n_0}\underline
{M}_\ell^{p} (\Pi_{n_1}u_1,\dots,
\Pi_{n_{p+1}}u_{p+1})=\Pi_{n_0}\widetilde{M}_\ell^{p}(\Pi_{n_1}u_1,\dots,
\Pi_{n_{p+1}}u_{p+1}),\end{align} where $F_m^\ell$ is defined by (\ref{eqn-G1}) and
$\widetilde{F}_m^{(e)}$ is defined by (\ref{eqn-G29}) with $e_0=\dots=e_\ell=e_{p+1}=-1,
e_{\ell+1}=\dots=e_{p}=1$.

When considering (\ref{eqn-H6}), we may assume n$_0\sim n_{p+1}$, $n_{p+1}\ge n'$ and
$(n_0,\dots, n_{p+1})\notin S_p^\ell$ if the right hand side of (\ref{eqn-H6}) is non
zero since we have (\ref{eqn-E1}) and (\ref{eqn-E13}) for $M_\ell^{p, \epsilon}$. Here
$S_p^\ell$ is the same as that in proposition \ref{thm-G6}. Thus the assumptions
concerning $(n_0, \dots, n_{p+1})$ in proposition \ref{thm-G6} hold true. We deduce from
(\ref{eqn-G32}) and the condition $\s{n_0}+\s{n_{p+1}}<\epsilon^{-\theta\kappa}$ that
\begin{equation}\label{eqn-H8}\begin{split}
&|F_m^\ell(\lambda_{n_0},\dots, \lambda_{n_{p+1}})|^{-1}\le
C(1+\s{n_0}+\s{n_{p+1}})^{3+\rho}
(1+\sqrt{n'})^{2N_0}\frac{S(n_0, \dots, n_{p+1})^{2N_0}}{\mu(n_0, \dots, n_{p+1})^{2N_0}} \\
&\le
C\epsilon^{-(4-a+\rho)\theta\kappa}(1+\sqrt{n_0}+\sqrt{n_{p+1}})^{a-1}(1+\sqrt{n'})^{2N_0}\frac{S(n_0,
\dots, n_{p+1})^{2N_0}}{\mu(n_0, \dots, n_{p+1})^{2N_0}}.
\end{split}\end{equation} for any $\rho>0$. Therefore if we define
\begin{equation}\label{eqn-H10}
\underline {M}_\ell^{p,\epsilon}(u_1,\dots,u_{p+1})=-\sum_{\substack{n\notin
S^\ell_p\\n_0\sim n_{p+1}, n_{p+1}\ge n'}}F_m^\ell(\lambda_{n_0},
\dots,\lambda_{n_{p+1}})^{-1}\Pi_{n_0}M_\ell^{p,\epsilon}(\Pi_{n_{1}}u_1, \dots,
\Pi_{n_{p+1}}u_{p+1}),
\end{equation} we obtain according to (\ref{eqn-H8}) and (\ref{eqn-E2}) that $\underline {M}_\ell^{p,\epsilon}\in
\widetilde{\mathcal{M}}_{p+1}^{\nu+\bar{\nu},2s-1}(\omega_\ell)$ with the first estimate
in (\ref{eqn-H4}) with $\bar{\nu}=2N_0$.

When considering (\ref{eqn-H6'}), we may assume $(n_0, \dots, n_{p+1})\notin N^{(e)}$
defined after (\ref{eqn-G29}). Actually, because of (\ref{eqn-E13}), we cannot find a
bijection $\sigma$ from $\{0,\dots, \ell, p+1\}$ to $\{\ell+1,\dots, p\}$ such that
$n_j=n_{\sigma(j)}$, $j=0,\dots,\ell,p+1$ if the right hand side of (\ref{eqn-H6'}) is
non zero. Consequently, we may use lower bound (\ref{eqn-G31}). If we define
$\underline{M}_\ell^p$ dividing in (\ref{eqn-H6'}) by $\widetilde{F}^{(e)}_m$, we thus
see that we get an element of $\underline
{M}_\ell^p\in\widetilde{\mathcal{M}}_{p+1}^{\nu+\bar{\nu},2s-2}(\widetilde{\omega}_\ell)$
for some $\bar{\nu}$. This completes the proof of (\ref{eqn-H3}) and (\ref{eqn-H4}).

(ii) We deduce again from (\ref{eqn-H5})
\begin{align}\label{eqn-H7'}
   -F_m^\ell(\lambda_{n_0}, \dots, \lambda_{n_{p+1}})\Pi_{n_0}\underline {R}_\ell^p(\Pi_{n_1}u_1,
   \dots, \Pi_{n_{p+1}}u_{p+1})=\Pi_{n_0}R
   _\ell^p  (\Pi_{n_1}u_1,\dots, \Pi_{n_{p+1}}u_{p+1}), \\\label{eqn-H8'}
   \widetilde{F}_m^{(e)}(\lambda_{n_0}, \dots, \lambda_{n_{p+1}})\Pi_{n_0}\underline
   {R}'^p_\ell(\Pi_{n_1}u_1,\dots,
   \Pi_{n_{p+1}}u_{p+1})=\Pi_{n_0}\widetilde{R}^p_\ell  (\Pi_{n_1}u_1,\dots,
   \Pi_{n_{p+1}}u_{p+1}),
\end{align} where $F_m^\ell$ and $\widetilde{F}_m^{(e)}$ are the same as in
(\ref{eqn-H6}) and (\ref{eqn-H6'}). Since $R_\ell^p\in
\widetilde{R}_{p+1}^{\nu,2s}(\omega_\ell)$ and thus (\ref{eqn-E21}) implies the right
hand side of (\ref{eqn-H7'}) vanishes if $(n_0, \dots, n_{p+1})\in S_p^\ell$, where
$S_p^\ell$ is defined in (\ref{eqn-F30}), we may assume $(n_0, \dots, n_{p+1})\notin
S^\ell_p$. Consequently, the condition of theorem \ref{thm-G1} is satisfied and we have
by (\ref{eqn-G3})
$$|F_m^\ell(\lambda_{n_0},\dots,\lambda_{n_{p+1}})|^{-1} \le
C(1+\s{n_0}+\s{n_1}+\dots+\s{n_{p+1}})^{2N_0+4}.$$  We then get an element of $\underline
{R}_\ell^p\in\widetilde{\mathcal{R}}_{p+1}^{\nu+\bar{\nu},2s}(\omega_\ell)$ dividing in
(\ref{eqn-H7'}) by $-F^\ell_m$ with $\bar{\nu}=2N_0+4$. Since
$\widetilde{R}_\ell^p\in\widetilde{R}_{p+1}^{\nu,2s}(\widetilde{\omega}_\ell)$, we see
that the right hand side of (\ref{eqn-H8'}) vanishes if $(n_0,\dots, n_{p+1})\in
\widetilde{S}^\ell_p$, where $\widetilde{S}^\ell_p$ is defined in (\ref{eqn-G2'}). This
implies that we may assume $(n_0, \dots, n_{p+1})\notin N^{(e)}$ which is defined after
(\ref{eqn-G29}) with $e_0=\dots=e_{\ell}=e_{p+1}=-1$, $e_{\ell+1}=\dots=e_p=1$. Thus the
condition of proposition \ref{thm-G5} is satisfied and we have
$$|\widetilde{F}_m^{(e)}(\lambda_{n_0}, \dots, \lambda_{n_{p+1}})|^{-1}\le
C(1+\sqrt{n_0}+\dots+\sqrt{n_{p+1}})^{N_0}.
$$ This allows us to get an element $\underline {R}'^p_\ell\in\widetilde
    {\mathcal{R}}_{p+1}^ {\nu+\bar{\nu},2s}(\widetilde{\omega}_\ell)$ for some $\bar{\nu}$
by dividing by $\widetilde{F}^{(e)}_m$ in (\ref{eqn-H8'}). This concludes the proof.
\end{proof}

\begin{pro}\label{thm-H3} Let $\mathcal {N}$ be the zero measure subset of $\R^*_+$ defined
in lemma \ref{thm-H1}, and fix $m\in\R^*_+-\mathcal{N}$. Let $\rho>0$ be any positive
number and $\Theta_s$ defined in (\ref{eqn-F4}).
 There are for any large enough integer $s$ , a map
$\Theta^1_s$, sending $\mathscr{H}^s(\mathbb{R}^d)\times (0, \frac12)$ to $\mathbb{R}$,
and maps $\ \Theta^2_s,\ \Theta^3_s,\ \Theta^4_s$ sending $\mathscr{H}^s(\mathbb{R}^d)$
to $\mathbb{R}$ such that there is a constant $C_s>0$ and for any $u\in
\mathscr{H}^s(\mathbb{R}^d)$ with $||u||_{\mathscr{H}^{s}}\le 1$ and any $\epsilon\in (0,
\frac12)$, one has
\begin{equation}\label{eqn-H12}\begin{split}
    &|\Theta^1_s(u, \epsilon)|\le C_s\epsilon^{-(4-a+\rho)\theta\kappa}||u||_{\mathscr{H}^
    {s}}^{\kappa+2},\qquad(a=2 \text{ if } d\ge 2
    \text{ and }\\&\qquad\qquad\qquad\qquad\qquad\qquad a=\frac{13}{6}-\varsigma \ \text
    { for any }\ \varsigma\in(0,1)\text{ if } d=1  ),\\
    &|\Theta^2_s(u)|, |\Theta_s^3(u)|,  |\Theta_s^4(u)|\le C_s||u||^{\kappa+2}_{\mathscr{H}^{s}}
\end{split}\end{equation} and such that
\begin{align}\label{eqn-H13}
        R(u)\overset{def}{=}\frac{d}{dt}\biggr[\Theta_s(u(t,\cdot))-\Theta^1_s(u(t,\cdot),
         \epsilon)-\Theta^2_s(u(t,\cdot))-\Theta^3_s(u(t,\cdot))-
        \Theta^4_s(u(t,\cdot))\biggr]
\end{align} satisfies
\begin{equation}\label{eqn-H14}
|R(u)|\le
C_s\epsilon^{-(4-a+\rho)\theta\kappa}||u||_{\mathscr{H}^s}^{2\kappa+2}+C_s\epsilon^{(a-1)\theta\kappa}
||u||^{\kappa+2}_{\mathscr{H}^{s}}+C_s||u||_{\mathscr{H}^{s}}^{2\kappa+2}.
\end{equation}
\end{pro}
\begin{proof} Considering the right hand side of (\ref{eqn-F5}), we decompose
\begin{equation}\label{eqn-H15}M^p_\ell(u_1, \dots, u_{p+1})=M^{p,\epsilon}_\ell(u_1,
\dots, u_{p+1})+V^{p,\epsilon}_\ell
(u_1, \dots, u_{p+1}),\end{equation} where the first term is given by (\ref{eqn-H2})
and the second one by
\begin{equation}\label{eqn-H16} V^{p, \epsilon}_\ell(u_1, \dots, u_{p+1})=\sum_{n_0}\sum_{n_{p+1}}\textbf{1}
_{\{\sqrt{n_0}+\sqrt{n_{p+1}}\ge \epsilon^{-\theta\kappa}\}}\Pi_{n_0}M^p_\ell(u_1, \dots, u_p, \Pi_{n_{p+1}}u_{p+1}).
\end{equation}
By definition \ref{def-E1}, we get for $a=2$ if $d\ge 2$ and $a=\frac{13}{6}-\varsigma$
if $d=1$
\begin{equation}\label{eqn-H17}\begin{split}||V^{p,\epsilon}_\ell(u_1, \dots, u_{p+1})||_{\mathscr {H}^
{-s}}\le C_N \sum_{n_0}\dots\sum_{p+1}(1+\sqrt{n_0}+\sqrt{n_{p+1}})^{2s-a}\frac{(1+\sqrt{n'})^{\nu}\mu
(n_0,\dots,n_{p+1})^N}{S(n_0, \dots, n_{p+1})^N}\\
\times\textbf{1}_{ \{\sqrt{n_0}+\sqrt{n_{p+1}}\ge \epsilon^{-\theta\kappa},
|n_{0}-n_{p+1}| <\frac12(n_{0}+n_{p+1}), n'\le
n_{p+1}\}}(1+\sqrt{n_0})^{-s}\prod_{j=1}^{p+1}||\Pi_{n_j}u_j||_{L^2}
\end{split}\end{equation} Following the proof of proposition \ref{thm-E1}, we know that
the gain of $a$ powers of
$\sqrt{n_0}+\sqrt{n_{p+1}}$ in the first term in the right hand side, coming from the
fact that $M^p_\ell\in \mathcal {M}^{\nu, 2s-a}_{p+1}$ , together with the condition
$\sqrt{n_0}+\sqrt{n_{p+1}}\ge \epsilon^{-\theta\kappa}$, allows us to estimate , for $N$
large enough and $s_0$ large enough with respect to $\nu$, (\ref{eqn-H17}) by
$C\epsilon^{(a-1)\theta\kappa}\Pi_{j=1}^p||u_j||_{\mathscr
{H}^{s_0}}||u_{p+1}||_{\mathscr {H}^s}$. Consequently, the quantity
\begin{equation}\label{eqn-H18} \sum_{p=\kappa}^{2\kappa-1}\sum_{\ell=0}^p Re\ i\langle V^{p,\epsilon}_\ell
(\bar{u},\dots,\bar{u}, u, \dots, u), u\rangle\end{equation} is bounded form above by the
second term of the right hand side of (\ref{eqn-H14}). In the rest of the proof, we may
therefore replace in the right hand side of (\ref{eqn-F5}) $M^p_\ell$ by
$M^{p,\epsilon}_\ell$.
\paragraph{}

Apply lemma \ref{thm-H1} to $M^{p,\epsilon}_\ell$, $\widetilde{M}^p_\ell$, $R^p_\ell$,
$\widetilde{R}^p_\ell$. This gives $\underline {M}^{p,\epsilon}_\ell, \underline
{M}^p_\ell, \underline {R}^p_\ell, \underline {R}'^p_\ell$. We set
\begin{equation}\label{eqn-H19}\begin{split}
\Theta^1_s(u(t,\cdot), \epsilon)=&\sum_{p=\kappa}^{2\kappa-1}\sum_{\ell=0}^{p}Re
\langle\underline {M}_\ell^{p,\epsilon}(\bar{u},\dots, \bar{u},u,\dots,u),u\rangle,\\
\Theta^2_s(u(t,\cdot))=&\sum_{p=\kappa}^{2\kappa-1}\sum_{\ell=0}^{p}Re \langle\underline
{M}_\ell^{p} (\bar{u},\dots, \bar{u},u,\dots,u, \bar{u}),u\rangle,\\
\Theta^3_s(u(t,\cdot))=&\sum_{p=\kappa}^{2\kappa-1}\sum_{\ell=0}^{p}Re \langle\underline
{R}_\ell^p(\bar{u},\dots,\bar{u},u,\dots,u),u\rangle,\\
\Theta^4_s(u(t,\cdot))=&\sum_{p=\kappa}^{2\kappa-1}\sum_{\ell=0}^{p}Re \langle\underline
{R}'^p_\ell(\bar{u},\dots,\bar{u},u,\dots,u, \bar{u}),u\rangle.
\end{split}\end{equation}
The general term in $\Theta^1_s(u(t,\cdot), \epsilon)$ has modulus bounded from above by
$$||\underline {M}_\ell^{p,\epsilon}(\bar{u},\dots, \bar{u},u,\dots, u)||_{\mathscr{H}^{-s}}
||u||_{\mathscr{H}^s}\le C\epsilon^{-
(4-a+\rho)\theta\kappa}||u||^\kappa_{\mathscr{H}^{s}}||u||^2_{\mathscr{H}^s}$$ for $u$ in
the unit ball of $\mathscr{H}^{s}(\mathbb{R}^d)$, using proposition \ref{thm-E1} with
$\tau=2s-1$ and proposition \ref{thm-B7} and (\ref{eqn-H4}). This gives the first
inequality of (\ref{eqn-H12}). To obtain the other estimates in (\ref{eqn-H12}), we apply
proposition \ref{thm-E1} to $\underline {M}^p_\ell$, remarking that if in (\ref{eqn-E3})
$\tau=2s-1$ and $s$ is large enough, the left hand side of (\ref{eqn-E3}) controls the
$\mathscr{H}^{-s}$ norm of $\underline {M}_\ell^p(\bar{u},\dots,\bar{u},u,\dots,u,
\bar{u})$. We also apply proposition \ref{thm-E2} with $\tau=2s$ in (\ref{eqn-E15}) to
$\underline {R}^p_\ell$, $\underline {R}'^p_\ell$. Then if $s_0$ is large enough, the
left hand side of (\ref{eqn-E15}) controls $\mathscr{H}^{-s}$ norm of $\underline
{R}^p_\ell(\bar{u}, \dots, \bar{u}, u, \dots, u)$ and  $\underline {R}'^p_\ell(\bar{u},
\dots, \bar{u}, u, \dots, u, \bar{u})$. These give us the other inequalities in
(\ref{eqn-H12}). Consequently we are left with proving (\ref{eqn-H14}). Remarking that we
may also write the equation as
\begin{equation}\label{eqn-H20} (D_t-\Lambda_m)u=-F\biggr(\Lambda_m^{-1}(\frac{u+\bar{u}}{2})\biggr),\end{equation}
we compute using notation (\ref{eqn-H1})
\begin{equation}\label{eqn-H21}\begin{split}
\frac{d}{dt}\Theta^1_s(u,\epsilon)=&\sum_{p=\kappa}^{2\kappa-1}\sum_{\ell=0}^{p}Re\
i\langle L^-_\ell(\underline {M}_\ell^{p,\epsilon})(\bar{u},\dots,\bar{u},u,\dots,u),
u\rangle\\
&+\sum_{p=\kappa}^{2\kappa-1}\sum_{\ell=0}^{p}\sum_{j=1}^\ell Re\ i\langle\underline
{M}_\ell^{p,\epsilon}(\bar{u},\dots,\bar{F},\dots,\bar{u},u,\dots,u),
u\rangle\\
&-\sum_{p=\kappa}^{2\kappa-1}\sum_{\ell=0}^{p}\sum_{j=\ell+1}^{p+1} Re\
i\langle\underline {M}_\ell^{p,\epsilon}(\bar{u},\dots,\bar{u},u,\dots,F,\dots,u),
u\rangle\\
&+\sum_{p=\kappa}^{2\kappa-1}\sum_{\ell=0}^{p} Re\ i\langle\underline
{M}_\ell^{p,\epsilon}(\bar{u},\dots,\bar{u},u,\dots,u), F\rangle.
\end{split}\end{equation}
By assumption on $F$, we have by proposition \ref{thm-B7} and \ref{thm-B9} that
$||F(v)||_{\mathscr{H}^s}\le C||u||_{\mathscr{H}^{s}}^\kappa||u||_{\mathscr{H}^s}$ if $s$
is large enough and $||u||_{\mathscr{H}^{s}}\le1$. Since $\underline {M}
_\ell^{p,\epsilon}\in\widetilde{\mathcal{M}}_{p+1}^{\nu+\bar{\nu},2s-1}(\omega_\ell)$, we
may apply proposition \ref{thm-E1} with $\tau=2s-1$ and (\ref{eqn-H4}) to see that the
last three terms in (\ref{eqn-H21}) have modulus bounded from above by the first term in
the right hand side of (\ref{eqn-H14}). When computing $\frac{d}{dt}\Theta_s(u)$, noting
that we have replaced $M^p_\ell$ by $M^{p,\epsilon}_\ell$, the first term in the right
hand side of (\ref{eqn-F5}) is the first term in the right hand side of (\ref{eqn-H21})
because of (\ref{eqn-H3}). Consequently, these contributions will cancel out each other
in the expression $\frac{d}{dt} [\Theta_s(u)-\Theta_s^1(u,\epsilon)]$. We compute
\begin{equation}\label{e14}\begin{split}
\frac{d}{dt}\Theta^2_s(u)=&\sum_{p=\kappa}^{2\kappa-1}\sum_{\ell=0}^{p}Re\ i\langle
L^+_\ell(\underline {M}_\ell^p)(\bar{u},\dots,\bar{u},u,
\dots,u,\bar{u}), u\rangle\\
&+\sum_{p=\kappa}^{2\kappa-1}\sum_{\ell=0}^{p}\sum_{j=1}^\ell Re\ i\langle\underline
{M}_\ell^p(\bar{u},\dots,\bar{F},\dots,
\bar{u},u,\dots,u, \bar{u}), u\rangle\\
&-\sum_{p=\kappa}^{2\kappa-1}\sum_{\ell=0}^{p}\sum_{j=\ell+1}^{p} Re\ i\langle\underline
{M}_\ell^p(\bar{u},\dots,\bar{u},u,\dots,F,\dots,u, \bar{u}),
u\rangle\\&+\sum_{p=\kappa}^{2\kappa-1}\sum_{\ell=0}^{p} Re\ i\langle\underline
{M}_\ell^p(\bar{u},,\dots,\bar{u},u,\dots,u,\bar{F}),
u\rangle\\
&+\sum_{p=\kappa}^{2\kappa-1}\sum_{\ell=0}^{p} Re\ i\langle\underline
{M}_\ell^p(\bar{u},,\dots,\bar{u},u,\dots,u,\bar{u}), F\rangle.
\end{split}\end{equation}
Since $\underline {M} _\ell^{p}\in\widetilde
{\mathcal{M}}_{p+1}^{\nu+\bar{\nu},2s-2}(\widetilde{\omega}_\ell)$, we have by
proposition \ref{thm-E1} with $\tau=2s-1$, proposition \ref{thm-B7} and (\ref{eqn-H4})
that the last three terms are estimated by the last term in the right hand side of
(\ref{eqn-H14}) if $s$ is large enough. The first one, according to lemma \ref{thm-H1},
cancels the contribution of $\widetilde{M}_\ell^p$ in (\ref{eqn-F5}) when computing
$R(u)$. We may treat $\Theta_s^3(u)$ and $\Theta_s^4(u)$ in the same way using
proposition \ref{thm-E2} with $\tau=2s$, and this will lead to the third term in the
right hand side of (\ref{eqn-H14}). Finally, the last term in (\ref{eqn-F5}) contributes
to the last term in the right hand side of (\ref{eqn-H14}). This concludes the proof of
the proposition.
\end{proof}
\paragraph{}

\begin{proof}[Proof of theorem 2.1.1]: We deduce from (\ref{eqn-H12}) and (\ref{eqn-H14})
\begin{gather}\label{e15}
\Theta_s(u(t,\cdot))\le \Theta_s(u(0,\cdot))-\Theta^1_s(u(0,\cdot),
\epsilon)-\Theta^2_s(u(0,\cdot))-\Theta^3_s(u(0,\cdot))-\Theta^4_s
(u(0,\cdot))\\\nonumber+\Theta^1_s(u(t,\cdot),
\epsilon)+\Theta^2_s(u(t,\cdot))+\Theta^3_s(u(t,\cdot))+\Theta^4_s(u(t,\cdot))\\\nonumber
+C_s\epsilon^{-(4-a+\rho)\theta\kappa}\int_0^t||u(t',\cdot)||_{\mathscr{H}^{s}}^{2
\kappa}||u(t',\cdot)||^2_{\mathscr{H}^s}dt'\\\nonumber
+C_s\epsilon^{(a-1)\theta\kappa}\int_0^t||u(t',\cdot)||_{\mathscr{H}^{s}}^{\kappa}
||u(t',\cdot)||^2_{\mathscr{H}^s}dt'\\\nonumber +C_s\int_0^t||u(t',
\cdot)||_{\mathscr{H}^{s}}^{2\kappa}||u(t', \cdot)||^2_{\mathscr{H}^s}dt',
\end{gather} where $a=2$  if  $d\ge 2$  and  $a=\frac{13}{6}-\varsigma$ for any
$\varsigma\in (0,1)$ if $d=1$. Take $\theta=\frac{1}{3+\rho}$ and $B>1$ a constant such
that for any $(v_0,v_1)$ in the unit ball of
$\mathscr{H}^{s+1}(\mathbb{R}^d)\times\mathscr{H}^s(\mathbb{R}^d)$,
$u(0,\cdot)=\epsilon(-iv_1+\Lambda_mv_0)$ satisfies $||u(0,\cdot)||_{\mathscr{H}^s}\le
B\epsilon$. Let $K>B$ be another constant to be chosen, and assume that for $\tau'$ in
some interval $[0, T]$ we have $||u(\tau',\cdot)||_{\mathscr{H}^s}\le K\epsilon\le 1$. If
$d\ge 2$, using (\ref{eqn-H12}) with $a=2$ we deduce from (\ref{e15}) and that there is a
constant $C>0$, independent of $B,K,\epsilon$, such that as long as $t\in[0, T]$
\begin{align*}
||u(t, \cdot)||^2_{\mathscr{H}^s}\le
C[B^2+\epsilon^{\frac{1}{3+\rho}\kappa}K^{\kappa+2}+t\epsilon^{\frac{4+\rho}{3+\rho}\kappa}(K^{2\kappa+2}
+K^{\kappa+2}) +t\epsilon^{2\kappa}K^{2\kappa+2}]\epsilon^2.
\end{align*}
If we assume that $T\le c\epsilon^{-\frac{4+\rho}{3+\rho}}$, where $\rho>0$ is arbitrary,
for a small enough $c>0$, and that $\epsilon$ is small enough, we get $||u(t,
\cdot)||_{\mathscr{H}^s}^2\le C(2B^2)\epsilon^2$. If $K$ has been chosen initially so
that $2CB^2<K^2$, we get by a standard continuity argument that the priori bound $||u(t,
\cdot)||_{\mathscr{H}^s}\le K\epsilon$ holds true on $[0,
c\epsilon^{-\frac{4+\rho}{3+\rho}}]$, in other words, the solution extends to such an
interval $|t|\le c\epsilon^{-\frac43(1-\rho)\kappa}$ with another arbitrary $\rho>0$. If
$d=1$, we may use (\ref{eqn-H12}) with $a=\frac{13}{6}-\varsigma$ to get
\begin{equation*} ||u(t, \cdot)||^2_{\mathscr{H}^s}\le
C[B^2+\epsilon^{\frac{7-6\varsigma}{18+6\rho}\kappa}K^{\kappa+2}+t\epsilon^{\frac{25+6
(\rho-\varsigma)}{18+6\rho}\kappa}(K^{2\kappa+2}+K^{\kappa+2})+t\epsilon^{2\kappa}K^{2
\kappa+2}]\epsilon^2.\end{equation*} With the same reasoning we may get in this case that
the solution extends to an interval of $|t|<c\epsilon^{-\frac{25}{18}(1-\rho)\kappa}$ for
some small $c>0$ and any $\rho>0$. This concludes the proof of the theorem.
\end{proof}

{\bf Acknowledgements}. The author thanks his advisors Daoyuan Fang and Jean-Marc Delort
for their guidance. Most of this work has been done during the stay of the author at
Universit\'{e} Paris-Nord, during the academic year 2007-2008.



\begin{thebibliography}{99}

\bibitem{Ba} D. Bambusi:  \emph{Birkhoff normal form for some nonlinear PDEs},
Comm. Math. Phys. 234 (2003), no. 2, 253--285.

\bibitem{BDGS} D. Bambusi, J.-M. Delort, B. Gr{\'e}bert and J. Szeftel: \emph{Almost global existence for Hamiltonian semi-linear
Klein-Gordon equations with small Cauchy data on Zoll manifolds}, Comm. Pure Appl. Math.
60 (2007), no. 11, 1665--1690.

\bibitem{BG} D. Bambusi and B. Gr{\'e}bert: \emph{Birkhoff normal form for
    partial differential equations with tame modulus},  Duke Math. J.
  135  (2006),  no. 3, 507--567.

\bibitem{BB} M. Berti and P. Bolle: \emph{Periodic solutions for higher dimensional nonlinear wave equations}, preprint (2007).

\bibitem{BM} E. Bierstone and P. Milman: \emph{Semianalytic and subanalytic sets}, Inst.
     Hautes, \'{E}tudes Sci. Publ. Math.  (1998), no. 67, 5-42.

\bibitem{Bo} J. Bourgain: \emph{Construction of approximative and
    almost periodic solutions of perturbed linear Schr{\"o}dinger and wave
    equations}, Geom. Funct. Anal. 6 (1996), no. 2, 201--230.

\bibitem{Bo2} J. Bourgain: \emph{Green's function estimates for lattice
    Schr{\"o}dinger operators and applications}. Annals of Mathematics
  Studies, 158. Princeton University Press, Princeton, NJ, (2005),
  x+173 pp.



\bibitem{Cr} W. Craig: \emph{Probl{\`e}mes de petits diviseurs dans les
    {\'e}quations aux d{\'e}riv{\'e}es partielles}.  Panoramas et Synth{\`e}ses,
  9. Soci{\'e}t{\'e}  Math{\'e}matique de France, Paris, (2000), viii+120 pp.

\bibitem{D} J.-M. Delort: \emph{On long time existence for small solutions of semi-linear\\
Klein-Gordon equations on the torus},  to appear, Journal Analyse Math{\'e}matique.

\bibitem{DS1} J.-M. Delort and J. Szeftel: \emph{Long-time existence for small data nonlinear Klein-Gordon equations on tori
        and spheres}, Internat. Math. Res. Notices (2004), no. 37, 1897-1966.

\bibitem{DS2} J.-M. Delort and J. Szeftel: \emph{Long-time existence for semi-linear Klein-Gordon equations with small Cauchy
        data on Zoll manifolds}, Amer. J. Math. 128 (2006), no. 5, 1187-1218.

\bibitem{EK} H. Eliasson and S. Kuksin: \emph{KAM For the non-linear
      Schr{\"o}dinger equation}, to appear, Annals of Mathematics.

\bibitem{G} B. Gr{\'e}bert: \emph{Birkhoff normal form and hamiltonian
    PDEs}, Partial differential equations and applications,  1--46, S\'emin. Congr., 15, Soc. Math. France, Paris, 2007.

\bibitem{GP} B. Gr{\'e}bert, R. Imekraz and \'E. Paturel: \emph{On the long time behavior for solutions of semi-linear harmonic
    oscillator with small Cauchy data on $R^d$}, preprint, (2008).

\bibitem{H} R. Hardt: \emph{Stratification of real analytic mappings and images}, Invent. Math. 28 (1975), 193-208.

\bibitem{He} B. Helffer: \emph{Th\'{e}orie spectrale pour des op\'{e}rateurs globalment
    elliptiques}, Ast\'{e}risque, 112  (1984).

\bibitem{Ho} L. H\"{o}rmander: \emph{Lectures on Nonlinear Hyperbolic Differential
    Equations}, Math\'{e}matiques \& Applications, Vol 26, Springer.


\bibitem{K} S. Klainerman: \emph{Global existence of small amplitude solutions to nonlinear
    Klein-Gordon equations in four space--time dimensions}, Comm. Pure Appl. Math. 38, (1985)
631--641.

\bibitem{KT} H. Koch and D. Tataru: \emph{$L^p$ eigenfunction bounds for the Hermite
    operator}, Duke Math. J. 128 (2005), no. 2, 369-392.


\bibitem{Ku1} S. Kuksin: \emph{Nearly integrable infinite-dimensional
    Hamiltonian systems}. Lecture Notes in Mathematics,
  1556. Springer-Verlag, Berlin, 1993. xxviii+101 pp.

\bibitem{Ku2} S. Kuksin: \emph{Analysis of Hamiltonian PDEs}. Oxford
  Lecture Series in Mathematics and its Applications, 19. Oxford
  University Press, Oxford, 2000. xii+212 pp.

\bibitem{OTT} T. Ozawa, K. Tsutaya and Y. Tsutsumi: \emph{Global
    existence and asymptotic behavior of solutions for the
    Klein-Gordon equations with quadratic nonlinearity in two space
    dimensions}.  Math. Z.  222  (1996),  no. 3, 341--362.

\bibitem{Sh} J. Shatah: \emph{Normal forms and quadratic nonlinear
    Klein-Gordon equations}.  Comm. Pure Appl. Math.  38  (1985),
  no. 5, 685--696.
\end{thebibliography}
\end{document}